# An essay on the general theory of stochastic processes*

## Ashkan Nikeghbali

*ETHZ*
*Departement Mathematik, Rämistrasse 101, HG G16*
*Zürich 8092, Switzerland*
*e-mail:* `ashkan.nikeghbali@math.ethz.ch`

**Abstract:** This text is a survey of the general theory of stochastic processes, with a view towards random times and enlargements of filtrations. The first five chapters present standard materials, which were developed by the French probability school and which are usually written in French. The material presented in the last three chapters is less standard and takes into account some recent developments.

## Contents

*This is an original survey paper

## 1. Introduction

P.A. Meyer and C. Dellacherie have created the so called general theory of stochastic processes, which consists of a number of fundamental operations on either real valued stochastic processes indexed by $[0, \infty)$, or random measures on $[0, \infty)$, relative to a given filtered probability space $\left(\Omega, \mathcal{F}, (\mathcal{F}_t)_{t \geq 0}, \mathbb{P}\right)$, where $(\mathcal{F}_t)$ is a right continuous filtration of $(\mathcal{F}, \mathbb{P})$ complete sub-$\sigma$-fields of $\mathcal{F}$.

This theory was gradually created from results which originated from the study of Markov processes, and martingales and additive functionals associated with them. A guiding principle for Meyer and Dellacherie was to understand to which extent the Markov property could be avoided; in fact, they were able to get rid of the Markov property in a radical way.

At this point, we would like to emphasize that, perhaps to the astonishment of some readers, stochastic calculus was not thought of as a basic "elementary" tool in 1972, when C. Dellacherie's little book appeared. Thus it seemed interesting to view some important facts of the general theory in relation with stochastic calculus.

The present essay falls into two parts: the first part, consisting of sections 2 to 5, is a review of the General Theory of Stochastic Processes and is fairly well known. The second part is a review of more recent results, and is much less so. Throughout this essay we try to illustrate as much as possible the results with examples.

More precisely, the plan of the essay is as follows:

- in Section 2, we recall the basic notions of the theory: stopping times, the optional and predictable $\sigma$-fields and processes,etc.
- in Section 3, we present the fundamental Section theorems;
- in Section 4, we present the fundamental Projection theorems;

- in Section 5, we recall the Doob-Meyer decomposition of semimartingales;
- in Section 6, we present a small theory of multiplicative decompositions of nonnegative local submartingales;
- in Section 7, we highlight the role of certain "hidden" martingales in the general theory of stochastic processes;
- in Section 8, we illustrate the theory with the study of arbitrary random times;
- in Section 9, we study how the basic operations depend on the underlying filtration, which leads us in fact to some introduction of the theory of enlargement of filtrations;

### Acknowledgements

I would like to thank an anonymous referee for his comments and suggestions which helped to improve the present text.

## 2. Basic notions of the general theory

**Throughout this essay, we assume we are given a filtered probability space $\left(\Omega, \mathcal{F}, (\mathcal{F}_t)_{t\geq 0}, \mathbb{P}\right)$ that satisfies the usual conditions**, that is $(\mathcal{F}_t)$ is a right continuous filtration of $(\mathcal{F}, \mathbb{P})$ complete sub-$\sigma$-fields of $\mathcal{F}$.

A stochastic process is said to be **càdlàg** if it almost surely has sample paths which are right continuous with left limits. A stochastic process is said to be **càglàd** if it almost surely has sample paths which are left continuous with right limits.

### 2.1. *Stopping times*

**Definition 2.1.** A stopping time is a mapping $T : \Omega \to \overline{\overline{\mathbb{R}}}_+$ such that $\{T \leq t\} \in \mathcal{F}_t$ for all $t \geq 0$.

To a given stopping time $T$, we associate the $\sigma$-field $\mathcal{F}_T$ defined by:

$$\mathcal{F}_T = \{A \in \mathcal{F},\ A \cap \{T \leq t\} \in \mathcal{F}_t \ for\ all\ t \geq 0\}.$$

We can also associate with $T$ the $\sigma$-field $\mathcal{F}_{T-}$ generated by $\mathcal{F}_0$ and sets of the form:

$$A \cap \{T > t\},\ with\ A \in \mathcal{F}_t\ and\ t \geq 0.$$

We recap here without proof some of the classical properties of stopping times.

**Proposition 2.2.** *Let $T$ be a stopping time. Then $T$ is measurable with respect to $\mathcal{F}_{T-}$ and $\mathcal{F}_{T-} \subset \mathcal{F}_T$.*

**Proposition 2.3.** *Let $T$ be a stopping time. If $A \in \mathcal{F}_T$, then*

$$T_A(\omega) = \begin{cases} T(\omega) & if\ \omega \in A \\ +\infty & if\ \omega \notin A \end{cases}$$

*is also a stopping time.*

**Proposition 2.4** ([26], Theorem 53, p.187). *Let $S$ and $T$ be two stopping times.*

1. *For every $A \in \mathcal{F}_S$, the set $A \cap \{S \leq T\} \in \mathcal{F}_T$.*
2. *For every $A \in \mathcal{F}_S$, the set $A \cap \{S < T\} \in \mathcal{F}_T-$.*

**Proposition 2.5** ([26], Theorem 56, p.189). *Let $S$ and $T$ be two stopping times such that $S \leq T$. Then $\mathcal{F}_S \subset \mathcal{F}_T$.*

One of the most used properties of stopping times is the optional stopping theorem.

**Theorem 2.6** ([69], Theorem 3.2, p.69). *Let $(M_t)$ be an $(\mathcal{F}_t)$ uniformly integrable martingale and let $T$ be a stopping time. Then, one has:*

$$\mathbb{E}\left[M_\infty \mid \mathcal{F}_T\right] = M_T \tag{2.1}$$

*and hence:*

$$\mathbb{E}\left[M_\infty\right] = \mathbb{E}\left[M_T\right] \tag{2.2}$$

One can naturally ask whether there exist some other random times (i.e. nonnegative random variables) such that (2.1) or (2.2) hold. We will answer these questions in subsequent sections.

### 2.2. Progressive, Optional and Predictable σ-fields

Now, we shall define the three fundamental $\sigma$-algebras we always deal with in the theory of stochastic processes.

**Definition 2.7.** A process $X = (X_t)_{t \geq 0}$ is called $(\mathcal{F}_t)$ progressive if for every $t \geq 0$, the restriction of $(t, \omega) \mapsto X_t(\omega)$ to $[0, t] \times \Omega$ is $\mathcal{B}[0, t] \otimes \mathcal{F}_t$ measurable. A set $A \in \mathbb{R}_+ \times \Omega$ is called progressive if the process $\mathbf{1}_A(t, \omega)$ is progressive. The set of all progressive sets is a $\sigma$-algebra called the progressive $\sigma$-algebra, which we will denote $\mathcal{M}$.

**Proposition 2.8** ([69], Proposition 4.9, p.44). *If $X$ is a $(\mathcal{F}_t)$ progressive process and $T$ is a $(\mathcal{F}_t)$ stopping time, then $X_T \mathbf{1}_{\{T < \infty\}}$ is $\mathcal{F}_T$ measurable.*

**Definition 2.9.** The optional $\sigma$-algebra $\mathcal{O}$ is the $\sigma$-algebra, defined on $\mathbb{R}_+ \times \Omega$, generated by all processes $(X_t)_{t \geq 0}$, adapted to $(\mathcal{F}_t)$, with càdlàg paths. A process $X = (X_t)_{t \geq 0}$ is called $(\mathcal{F}_t)$ optional if the map $(t, \omega) \mapsto X_t(\omega)$ is measurable with respect to the optional $\sigma$-algebra $\mathcal{O}$.

**Proposition 2.10.** *Let $X = (X_t)_{t \geq 0}$ be an optional process and $T$ a stopping time. Then:*

1. *$X_T \mathbf{1}_{\{T < \infty\}}$ is $\mathcal{F}_T$ measurable.*
2. *the stopped process $X^T = (X_{t \wedge T})_{t \geq 0}$ is optional.*

**Definition 2.11.** The predictable $\sigma$-algebra $\mathcal{P}$ is the $\sigma$-algebra, defined on $\mathbb{R}_+ \times \Omega$, generated by all processes $(X_t)_{t \geq 0}$, adapted to $(\mathcal{F}_t)$, with left continuous paths on $]0, \infty[$. A process $X = (X_t)_{t \geq 0}$ is called $(\mathcal{F}_t)$ predictable if the map $(t, \omega) \mapsto X_t(\omega)$ is measurable with respect to the predictable $\sigma$-algebra $\mathcal{P}$.

**Proposition 2.12** ([26], Theorem 67, p.200). *Let $X = (X_t)_{t \geq 0}$ be a predictable process and $T$ a stopping time. Then:*

1. $X_T \mathbf{1}_{\{T < \infty\}}$ *is $\mathcal{F}_T-$ measurable.*
2. *the stopped process $X^T = (X_{t \wedge T})_{t \geq 0}$ is predictable.*

The following inclusions always hold:

$$\mathcal{P} \subset \mathcal{O} \subset \mathcal{M}.$$

It is easy to show that every adapted càdlàg process is progressively measurable ([69]), hence $\mathcal{O} \subset \mathcal{M}$. We also have that every càg and adapted process is predictable (every such process can be written as a limit of càdlàg processes), thus $\mathcal{P} \subset \mathcal{O}$.

The inclusions may be strict.

**Example 2.13.** $\mathcal{O} \subsetneq \mathcal{M}$. The following example is due Dellacherie and Meyer (see [24], p.128). Let $B = (B_t)_{t \geq 0}$ be the standard Brownian motion and let $(\mathcal{F}_t)$ be the usual augmentation of the natural filtration of $B$. For each $\omega$, the set

$$\{s : \quad B_s(\omega) \neq 0\}$$

is the disjoint union of open intervals, the excursions intervals. Define the set $E$ by:

$$E \equiv \{(s, \omega) : \quad s \text{ is the left-hand endpoint of an excursion interval}\}.$$

Then $E$ is progressive but not optional.

**Example 2.14.** $\mathcal{P} \subsetneq \mathcal{O}$. The standard Poisson process $N = (N_t)_{t \geq 0}$ is optional, but not predictable, in its own filtration $(\mathcal{F}_t)$.

Now, we give a characterization of the optional and predictable $\sigma$-algebras in terms of stochastic intervals.

**Definition 2.15.** Let $S$ and $T$ be two positive random variables such that $S \leq T$. We define $[S, T[$, as the random subset of $\mathbb{R}_+ \times \Omega$:

$$[S, T[ = \{(t, \omega) : \quad S(\omega) \leq t < T(\omega)\}.$$

All the stochastic intervals we can form with a pair of stopping times $(S, T)$, such that $S \leq T$, are optional. Indeed, $\mathbf{1}_{[0,T[}$ and $\mathbf{1}_{[S,\infty[}$ are càdlàg. and adapted: hence $\mathbf{1}_{[S,T[}$ is optional. Taking $T = S + \frac{1}{n}$ and letting $n$ go to infinity yields $[S]$ is optional. It is therefore immediate that the other type of stochastic intervals are optional.

Now let $Z$ be an $\mathcal{F}_S$ measurable random variable; then the process $Z(\omega) \mathbf{1}_{[S,T]}(t,\omega)$ is optional since it is adapted and càdlàg. We can also prove an analogous result with other types of stochastic intervals.

The same way as above, it can be shown that $]S,T]$ is predictable and that for $Z$ an $\mathcal{F}_S$ measurable random variable, $Z(\omega) \mathbf{1}_{]S,T]}(t,\omega)$ is predictable.

Now, we can give a characterization of the optional and predictable $\sigma$-algebras in terms of stochastic intervals.

**Proposition 2.16.** *The optional $\sigma$-algebra is generated by stochastic intervals of the form $[T,\infty[$, where $T$ is a stopping time:*

$$\mathcal{O} = \sigma\{[T,\infty[, \quad T \text{ is a stopping time}\}.$$

*Moreover, if $Y$ is a $\mathcal{F}_T$ measurable random variable, then there exists an optional process $(X_t)_{t \in \overline{\mathbb{R}}_+}$ such that $Y = X_T$.*

*Proof.* From the remark above, it suffices to show that any càdlàg and adapted process $X$ can be approximated with a sequence of elements of $\sigma\{[T,\infty[,\ T$ is a stopping time$\}$. Fix $\varepsilon > 0$ and let: $T_0 \equiv 0$ and $Z_0 = X_{T_0}$. Now define inductively for $n \geq 0$,

$$\begin{aligned} T_{n+1} &= \inf\{t > T_n : |X_t - Z_n| > \varepsilon\} \\ Z_{n+1} &= X_{T_{n+1}} \quad \text{if } T_{n+1} < \infty. \end{aligned}$$

Since $X$ has left limits, $T_n \uparrow \infty$. Now set:

$$Y \equiv \sum_{n \geq 0} Z_n \mathbf{1}_{[T_n, T_{n+1}[}.$$

Then $|X - Y| \leq \varepsilon$ and this completes the proof. $\qquad\square$

**Remark 2.17.** We also have:

$$\mathcal{O} = \sigma\{[0,T[, \quad T \text{ is a stopping time}\}.$$

**Remark 2.18.** It is useful to note that for a random time $T$, we have that $[T,\infty[$ is in the optional sigma field if and only if $T$ is a stopping time.

A similar result holds for the predictable $\sigma$-algebra (see [26], [39] or [69]).

**Proposition 2.19** ([26], Theorem 67, p. 200)**.** *The predictable $\sigma$-algebra is generated by one of the following collections of random sets:*

1. *$A \times \{0\}$ where $A \in \mathcal{F}_0$, and $[0,T]$ where $T$ is a stopping time;*
2. *$A \times \{0\}$ where $A \in \mathcal{F}_0$, and $A \times (s,t]$ where $s < t$, $A \in \mathcal{F}_s$;*

Now we give an easy result which is often used in martingale theory.

**Proposition 2.20.** *Let $X = (X_t)_{t \geq 0}$ be an optional process. Then:*

1. *The jump process $\Delta X \equiv X - X_-$ is optional;*
2. *$X_-$ is predictable;*
3. *if moreover $X$ is predictable, then $\Delta X$ is predictable.*

### 2.3. Classification of stopping times

We shall now give finer results about stopping times. The notions developed here are very useful in the study of discontinuous semimartingales (see [39] for example). The proofs of the results presented here can be found in [24] or [26].

We first introduce the concept of predictable stopping times.

**Definition 2.21.** A predictable time is a mapping $T : \Omega \to \overline{\mathbb{R}}_+$ such that the stochastic interval $[0, T[$ is predictable.

Every predictable time is a stopping time since $[T, \infty[ \in \mathcal{P} \subset \mathcal{O}$. Moreover, as $[T] = [0, T] \setminus [0, T[$, we deduce that $[T] \in \mathcal{P}$.

We also have the following characterization of predictable times:

**Proposition 2.22** ([26], Theorem 71, p.204)**.** *A stopping time $T$ is predictable if there exists a sequence of stopping times $(T_n)$ satisfying the following conditions:*

  *1. $(T_n)$ is increasing with limit $T$.*
  *2. we have $\{T_n < T\}$ for all $n$ on the set $\{T > 0\}$;*

  *The sequence $(T_n)$ is called an announcing sequence for $T$.*

Now we enumerate some important properties of predictable stopping times, which can be found in [24] p.54, or [26] p.205.

**Theorem 2.23.** *Let $S$ be a predictable stopping time and $T$ any stopping time. For all $A \in \mathcal{F}_{S-}$, the set $A \cap \{S \leq T\} \in \mathcal{F}_{T-}$. In particular, the sets $\{S \leq T\}$ and $\{S = T\}$ are in $\mathcal{F}_{T-}$.*

**Proposition 2.24.** *Let $S$ and $T$ be two predictable stopping times. Then the stopping times $S \wedge T$ and $S \vee T$ are also predictable.*

**Proposition 2.25.** *Let $A \in \mathcal{F}_{T-}$ and $T$ a predictable stopping time. Then the time $T_A$ is also predictable.*

**Proposition 2.26.** *Let $(T_n)$ be an increasing sequence of predictable stopping times and $T = \lim_n T_n$. Then $T$ is predictable.*

We recall that a random set $A$ is called evanescent if the set

$$\{\omega : \ \exists \ t \in \mathbb{R}_+ \text{ with } \ (t, \omega) \in A\}$$

is $\mathbb{P}-$null.

**Definition 2.27.** Let $T$ be a stopping time.

1. We say that $T$ is accessible if there exists a sequence $(T_n)$ of predictable stopping times such that:

$$[T] \subset (\cup_n [T_n])$$

up to an evanescent set, or in other words,

$$\mathbb{P} \left[ \cup_n \{\omega : \ T_n(\omega) = T(\omega) < \infty\} \right] = 1$$

2. We say that $T$ is totally inaccessible if for all predictable stopping times $S$ we have:

$$[T] \cap [S] = \varnothing$$

up to an evanescent set, or in other words:

$$\mathbb{P}\left[\{\omega : \ T(\omega) = S(\omega) < \infty\}\right] = 0.$$

**Remark 2.28.** It is obvious that predictable stopping times are accessible and that the stopping times which are both accessible and totally inaccessible are almost surely infinite.

**Remark 2.29.** There exist stopping times which are accessible but not predictable.

**Theorem 2.30** ([26] Theorem 81, p.215). *Let $T$ be a stopping time. There exists a unique (up to a $\mathbb{P}-null$ set) partition of the set $\{T < \infty\}$ into two sets $A$ and $B$ which belong to $\mathcal{F}_{T-}$ such that $T_A$ is accessible and $T_B$ is totally inaccessible. The stopping time $T_A$ is called the accessible part of $T$ while $T_B$ is called the totally inaccessible part of $T$.*

Now let us examine a special case where the accessible times are predictable. For this, we need to define the concept of quasi-left continuous filtrations.

**Definition 2.31.** The filtration $(\mathcal{F}_t)$ is quasi-left continuous if

$$\mathcal{F}_T = \mathcal{F}_{T-}$$

for all predictable stopping times.

**Theorem 2.32** ([26] Theorem 83, p.217). *The following assertions are equivalent:*

1. *The accessible stopping times are predictable;*
2. *The filtration $(\mathcal{F}_t)$ is quasi-left continuous;*
3. *The filtration $(\mathcal{F}_t)$ does not have any discontinuity time:*

$$\bigvee \mathcal{F}_{T_n} = \mathcal{F}_{(\lim T_n)}$$

*for all increasing sequences of stopping times $(T_n)$.*

**Definition 2.33.** A càdlàg process $X$ is called quasi-left continuous if $\Delta X_T = 0$, a.s. on the set $\{T < \infty\}$ for every predictable time $T$.

**Definition 2.34.** A random set $A$ is called thin if it is of the form $A = \cup [T_n]$, where $(T_n)$ is a sequence of stopping times; if moreover the sequence $(T_n)$ satisfies $[T_n] \cap [T_m] = \varnothing$ for all $n \neq m$, it is called an exhausting sequence for $A$.

**Proposition 2.35.** *Let $X$ be a càdlàg adapted process. The following are equivalent:*

1. *$X$ is quasi-left continuous;*

2. *there exists a sequence of totally inaccessible stopping times that exhausts the jumps of $X$;*
3. *for any increasing sequence of stopping times $(T_n)$ with limit $T$, we have $\lim X_{T_n} = X_T$ a.s. on the set $\{T < \infty\}$.*

### *2.4. Début theorems*

In this section, we give a fundamental result for realizations of stopping times: the début theorem. Its proof is difficult and uses the same hard theory (capacities theory) as the section theorems which we shall state in the next section.

**Definition 2.36.** Let $A$ be a subset of $\mathbb{R}_+ \times \Omega$. The début of $A$ is the function $D_A$ defined as:

$$D_A(\omega) = \inf \{t \in \mathbb{R}_+ : \ (t, \omega) \in A\},$$

with $D_A(\omega) = \infty$ if this set is empty.

It is a nice and difficult result that when the set $A$ is progressive, then $D_A$ is a stopping time ([24], [26]):

**Theorem 2.37** ([24], Theorem 23, p. 51)**.** *Let $A$ be a progressive set, then $D_A$ is a stopping time.*

Conversely, every stopping time is the début of a progressive (in fact optional) set: indeed, it suffices to take $A = [T, \infty[$ or $A = [T]$.

The proof of the début theorem is an easy consequence of the following difficult result from measure theory:

**Theorem 2.38.** *If $(E, \mathcal{E})$ is a locally compact space with a countable basis with its Borel $\sigma$-field and $(\Omega, \mathcal{F}, \mathbb{P})$ is a complete probability space, for every set $A \in \mathcal{E} \otimes \mathcal{F}$, the projection $\pi(A)$ of $A$ into $\Omega$ belongs to $\mathcal{F}$.*

*Proof of the début theorem.* We apply Theorem 2.38 to the set $A_t = A \cap ([0, t[ \times \Omega)$ which belongs to $\mathcal{B}([0, t[) \otimes \mathcal{F}_t$. As a result, $\{D_A \leq t\} = \pi(A_t)$ belongs to $\mathcal{F}_t$. $\square$

We can define the n-début of a set $A$ by

$$D_A^n(\omega) = \inf \{t \in \mathbb{R}_+ : \ [0, t] \cap A \text{ contains at least n points}\};$$

we can also define the $\infty$-début of $A$ by:

$$D_A^n(\omega) = \inf \{t \in \mathbb{R}_+ : \ [0, t] \cap A \text{ contains infinitely many points}\}.$$

**Theorem 2.39.** *The n-début of a progressive set $A$ is a stopping time for $n = 1, 2, \ldots, \infty$.*

*Proof.* The proof is easy once we know that $D_A^1(\omega)$ is a stopping time. Indeed, by induction on $n$, we prove that $D_A^{n+1}(\omega)$ which is the début of the progressive set $A_n = A \cap ]D_A^n(\omega), \infty[$. $D_A^\infty(\omega)$ is also a stopping time as the début of the progressive set $\cap A_n$. $\square$

It is also possible to show that the penetration time $T$ of a progressive set $A$, defined by:

$$T(\omega) = \inf\left\{t \in \mathbb{R}_+ : \ [0,t] \cap A \text{ contains infinitely non countable many points}\right\}$$

is a stopping time.

We can naturally wonder if the début of a predictable set is a predictable stopping time. One moment of reflexion shows that the answer is negative: every stopping time is the début of the predictable set $]T, \infty[$ without being predictable itself. However, we have:

**Proposition 2.40.** *Let $D_A$ be the début of a predictable set $A$. If $[D_A] \subset A$, then $D_A$ is a predictable stopping time.*

*Proof.* If $[D_A] \subset A$, then $[D_A] = A \cap [0, D_A]$, is predictable since $A$ is predictable and $D_A$ is a stopping time. Hence $D_A$ is predictable. $\square$

One can deduce from there that:

**Proposition 2.41.** *Let $A$ be a predictable set which is closed for the right topology[1]. Then its début $D_A$ is a predictable stopping time.*

Now we are going to link the above mentioned notions to the jumps of some stochastic processes. We will follow [39], chapter I.

**Lemma 2.42.** *Any thin random set admits an exhausting sequence of stopping times.*

**Proposition 2.43.** *If $X$ is a càdlàg adapted process, the random set $U \equiv \{\Delta X \neq 0\}$ is thin; an exhausting sequence $(T_n)$ for this set is called a sequence that exhausts the jumps of $X$. Moreover, if $X$ is predictable, the stopping times $(T_n)$ can be chosen predictable.*

*Proof.* Let
$$U_n \equiv \left\{(t, \omega) : \ |X_t(\omega) - X_{t-}(\omega) > 2^{-n}\right\},$$

for $n$ an integer and set $V_0 = U_0$ and

$$V_n = U_n - U_{n-1}.$$

The sets $V_n$ are optional (resp. predictable if $X$ is predictable) and are disjoint. Now, let us define the stopping times

$$
\begin{aligned}
D_n^1 &= \inf\left\{t : \ (t, \omega) \in V_n\right\} \\
D_n^{k+1} &= \inf\left\{t > D_n^k : \ (t, \omega) \in V_n\right\},
\end{aligned}
$$

[1]We recall that the right topology on the real line is a topology whose basis is given by intervals $[s, t[$.

so that $D_n^j$ represents the $j$—th jump of $X$ whose size in absolute value is between $2^{-n}$ and $2^{-n+1}$. Since $X$ is càdlàg, $V_n$ does not have any accumulation point and the stopping times $(D_n^k)_{(k,n) \in \mathbb{N}^2}$ enumerate all the points in $V_n$. Moreover, from Proposition 2.41, the stopping times $(D_n^k)$ are predictable if $X$ is predictable. To complete the proof, it suffices to index the doubly indexed sequence $(D_n^k)$ into a simple indexed sequence $(T_n)$. $\qquad \square$

In fact, we have the following characterization for predictable processes:

**Proposition 2.44.** *If $X$ is càdlàg adapted process, then $X$ is predictable if and only if the following two conditions are satisfied:*

1. *For all totally inaccessible stopping times $T$,*

$$\Delta X_T = 0, \text{ a.s.on} \{T < \infty\}$$

2. *For every predictable stopping time $T$, $X_T \mathbf{1}_{\{T < \infty\}}$ is $\mathcal{F}_{T-}$ measurable.*

Finally, we characterize $\mathcal{F}_T$ and $\mathcal{F}_{T-}$ measurable random variables:

**Theorem 2.45.** *Let $T$ be a stopping time. A random variable $Z$ is $\mathcal{F}_{T-}$ measurable if and only if there exists a predictable process $(X_t)$ such that $Z = X_T$ on the set $\{T < \infty\}$. Similarly, a random variable $Z$ is $\mathcal{F}_T$ measurable if and only if there exists an optional process $(Y_t)$ such that $Z = Y_T$ on the set $\{T < \infty\}$.*

*Proof.* We only prove the first part, the proof for the second part being the same. We have just seen that the condition is sufficient. To show the condition is necessary, by the monotone class theorem, it suffices to check the theorem for indicators of sets that generate the sigma field $\mathcal{F}_{T-}$, i.e. for $Z = \mathbf{1}_A$ when $A \in \mathcal{F}_0$ and for $Z = \mathbf{1}_{B \cap \{s < T\}}$, where $B \in \mathcal{F}_s$. But then, one can take $X = \mathbf{1}_{[0_A, \infty[}$ in the first case and $X = \mathbf{1}_{]s_B, \infty[}$ in the second case. $\qquad \square$

## 3. Section theorems

This section is devoted to a deep and very difficult result called the *section theorem*. The reader can refer to [26], p.219-220 or [24], p.70 for a proof. We will illustrate the theorem with some standard examples (here again the examples we deal with can be found in [26], [24] or [70]).

**Theorem 3.1** (Optional and predictable section theorems). *Let $A$ be an optional (resp. predictable) set. For every $\varepsilon > 0$, there is a stopping time (resp. predictable stopping time) $T$ such that:*

1. *$[T] \subset A$,*
2. *$\mathbb{P}[T < \infty] \geq \mathbb{P}(\pi(A)) - \varepsilon$,*

*where $\pi$ is the canonical projection of $\Omega \times \mathbb{R}_+$ onto $\Omega$.*

Throughout the paper, we shall use the optional and predictable section theorems. For now, we give some classical applications.

**Theorem 3.2.** *Let $(X_t)$ and $(Y_t)$ be two optional (resp. predictable) processes. If for every finite stopping time (resp. every finite predictable stopping time) one has:*

$$X_T = Y_T \quad \text{a.s.,}$$

*then the processes $(X_t)$ and $(Y_t)$ are indistinguishable.*

*Proof.* It suffices to apply the section theorem to the optional (resp. predictable) set

$$A = \{(t,\omega): \ X_t(\omega) \neq Y_t(\omega)\}.$$

Indeed, if the set $A$ were not evanescent, there would exist a stopping time whose graph would not be evanescent and which would be contained in $A$. This would imply the existence of some $t \in \mathbb{R}_+$ such that $X_{T \wedge t}$ would not be equal to $Y_{T \wedge t}$ almost surely. $\qquad\square$

**Theorem 3.3.** *Let $(X_t)$ and $(Y_t)$ be two optional (resp. predictable) processes. If for every stopping time (resp. every predictable stopping time) one has:*

$$\mathbb{E}\left[X_T \mathbf{1}_{T<\infty}\right] = \mathbb{E}\left[Y_T \mathbf{1}_{T<\infty}\right], \tag{3.1}$$

*then the processes $(X_t)$ and $(Y_t)$ are indistinguishable[2]*

*Proof.* It suffices to apply the section theorems to the sets:

$$\begin{aligned}
A &= \{(t,\omega): \ X_t(\omega) < Y_t(\omega)\} \\
A' &= \{(t,\omega): \ X_t(\omega) > Y_t(\omega)\}.
\end{aligned}$$

$\square$

**Remark 3.4.** It is very important to check (3.1) for any stopping time, finite or not. Indeed, if $(M_t)$ is a uniformly integrable martingale, with $M_0 = 0$, then for every finite stopping time,

$$\mathbb{E}\left[M_T \mathbf{1}_{T<\infty}\right] = 0,$$

but $(M_t)$ is not usually indistinguishable from the null process.

To conclude this section, we give two other well known results as a consequence of the section theorems. We first recall the definition of the class $(D)$:

**Definition 3.5** (class $(D)$)**.** A process $X$ is said to be of class $(D)$ if the family

$$\{X_T \mathbf{1}_{T<\infty}, \ T \text{ a stopping time}\}$$

is uniformly integrable ($T$ ranges through all stopping times).

---

[2] Two processes $X$ and $X'$ defined on the same probability space are called indistinguishable if for almost all $\omega$,

$$X_t = X'_t \qquad \text{for every t.}$$

**Proposition 3.6.** *Let* $(Z_t)_{0 \le t \le \infty}$ *be an optional process. Assume that for all stopping times* $T$, *the random variable* $Z_T$ *is in* $L^1$ *and* $\mathbb{E}[Z_T]$ *does not depend on* $T$. *Then* $(Z_t)$ *is a uniformly integrable martingale which is right continuous (up to an evanescent set).*

*Proof.* Let $T$ be a stopping time and let $A \in \mathcal{F}_T$. The assumptions of the proposition yield:

$$\mathbb{E}[Z_{T_A}] = \mathbb{E}[Z_\infty],$$

and hence

$$\int_A Z_T d\mathbb{P} + \int_{A^c} Z_\infty d\mathbb{P} = \int_A Z_\infty d\mathbb{P} + \int_{A^c} Z_\infty d\mathbb{P}.$$

Consequently, we have:

$$Z_T = \mathbb{E}[Z_\infty | \mathcal{F}_T] \text{ a.s.}$$

Now define also $(X_t)$ as the càdlàg version of the martingale:

$$X_t = \mathbb{E}[Z_\infty | \mathcal{F}_t].$$

By the optional stopping theorem, we have

$$X_T = \mathbb{E}[Z_\infty | \mathcal{F}_T],$$

and with an application of the section theorem, we obtain that $X$ and $Z$ are indistinguishable and this completes the proof of the proposition. □

## 4. Projection theorems

In this section we introduce the fundamental notions of optional (resp. predictable) projection and dual optional (resp. predictable) projection. These projections play a very important role in the general theory of stochastic processes. We shall give some nice applications in subsequent sections (in particular we shall see how the knowledge of the dual predictable projection of some honest times may lead to quick proofs of multidimensional extensions of Paul Lévy's arc sine law).

Here again, the material is standard in the general theory of stochastic processes and the reader can refer to the books [24] or [26, 27] for more details and refinements.

### *4.1. The optional and predictable projections*

By convention, we take $\mathcal{F}_0 = \mathcal{F}_{0-}$.

**Theorem 4.1.** *Let* $X$ *be a measurable process either positive or bounded. There exists a unique (up to indistinguishability) optional process* $Y$ *such that:*

$$\mathbb{E}[X_T \mathbf{1}_{T<\infty} | \mathcal{F}_T] = Y_T \mathbf{1}_{T<\infty} \text{ a.s.}$$

*for every stopping time* $T$.

**Definition 4.2.** The process $Y$ is called the *optional projection* of $X$, and is denoted by $^oX$.

*Proof.* The uniqueness follows from the optional section theorem. The space of bounded processes $X$ which admit an optional projection is a vector space. Moreover, let $X^n$ be a uniformly bounded increasing sequence of processes with limit $X$ and suppose they admit projections $Y^n$. The section theorem again shows that the sequence $(Y^n)$ is a.s. increasing; it is easily checked that $\varliminf Y^n$ is a projection for $X$.

By the monotone class theorem (as it is stated for example in [69] Theorem 2.2, p.3) it is now enough to prove the statement for a class of processes closed under pointwise multiplication and generating the sigma field $\mathcal{F} \otimes \mathcal{B}(\mathbb{R}_+)$. Such a class is provided by the processes

$$X_t(\omega) = \mathbf{1}_{[0,u[}(t)H(\omega), \quad 0 \leq u \leq \infty, \ H \in L^\infty(\mathcal{F}).$$

Let $(H_t)$ be a càdlàg version of $\mathbb{E}[H|\mathcal{F}_t]$ (with the convention that $H_{0-} = H_0$). The optional stopping theorem proves that

$$Y_t = \mathbf{1}_{[0,u[}(t)H_t$$

satisfies the condition of the statement. The proof is complete in the bounded case. For the general case we use the processes $X \wedge n$ and pass to the limit. $\quad\square$

**Theorem 4.3.** *Let $X$ be a measurable process either positive or bounded. There exists a unique (up to indistinguishability) predictable process $Z$ such that:*

$$\mathbb{E}[X_T\mathbf{1}_{T<\infty}|\mathcal{F}_{T-}] = Z_T\mathbf{1}_{T<\infty} \text{ a.s.}$$

*for every predictable stopping time $T$.*

**Definition 4.4.** The process $Z$ is called the *predictable projection* of $X$, and is denoted by $^pX$.

*Proof.* The proof is exactly the same as the proof of Theorem 4.1, except that at the end, one has to apply the predictable stopping theorem we shall give below. $\quad\square$

**Theorem 4.5** (Predictable stopping theorem). *Let $(X_t)$ be a right continuous and uniformly integrable martingale. Then for any predictable stopping time, we have:*

$$X_{T-} = \mathbb{E}[X_T|\mathcal{F}_{T-}] = \mathbb{E}[X_\infty|\mathcal{F}_{T-}].$$

*Consequently, the predictable projection of $X$ is the process $(X_{t-})$.*

*Proof.* Let $(T_n)$ be a sequence of stopping times that announces $T$; we then have

$$\begin{aligned} X_{T-} &= \lim_{n \to \infty} X_{T_n} = \mathbb{E}\left[X_T | \bigvee_n \mathcal{F}_{T_n}\right] \\ &= \mathbb{E}\left[X_\infty | \bigvee_n \mathcal{F}_{T_n}\right]. \end{aligned}$$

Since $\mathcal{F}_{T-} = \bigvee_n \mathcal{F}_{T_n}$, the desired result follows easily. $\qquad\blacksquare$

**Remark 4.6.** In the proof of Theorem 4.1, we have in fact given some fundamental examples of projections: if $H$ is an integrable random variable, and if $(H_t)$ is a càdlàg version of the martingale $\mathbb{E}[H|\mathcal{F}_t]$, then the optional projection of the constant process $X_t(\omega) = H(\omega)$ is $(H_t)$ and its predictable projection is $(H_{t-})$.

**Corollary 4.7.** *If $X$ is a predictable local martingale, then $X$ is continuous.*

*Proof.* By localization, it suffices to consider uniformly integrable martingales. Let $T$ be a predictable stopping time. Since $X$ is predictable, $X_T$ is $\mathcal{F}_{T-}$ measurable and hence: $\mathbb{E}[X_T|\mathcal{F}_{T-}] = X_T$, and from the predictable stopping theorem, $X_T = X_{T-}$. As this holds for any predictable stopping time and $X$ is predictable, we conclude that $X$ has no jumps. $\qquad\blacksquare$

We can also mention another corollary of the predictable stopping theorem which is not so well known:

**Corollary 4.8** ([27], p. 99)**.** *Let $(X_t)$ be a right continuous local martingale. Then for any predictable stopping time, we have:*

$$\mathbb{E}[|X_T|\mathbf{1}_{T<\infty}|\mathcal{F}_{T-}] < \infty \text{ a.s.,}$$

*and*

$$\mathbb{E}[X_T\mathbf{1}_{T<\infty}|\mathcal{F}_{T-}] = X_{T-}\,\mathbf{1}_{T<\infty}.$$

By considering $T_\Lambda$, where $\Lambda \in \mathcal{F}_T$ or $\Lambda \in \mathcal{F}_{T-}$, and where $T$ is a stopping time or a predictable stopping, we may as well restate Theorems 4.1 and 4.3 as:

**Theorem 4.9.** *Let $X$ be a measurable process either positive or bounded.*

1. *There exists a unique (up to indistinguishability) optional process $Y$ such that:*
$$\mathbb{E}[X_T\mathbf{1}_{T<\infty}] = \mathbb{E}[Y_T\mathbf{1}_{T<\infty}]$$
   *or any stopping time $T$.*
2. *There exists a unique (up to indistinguishability) predictable process $Z$ such that:*
$$\mathbb{E}[X_T\mathbf{1}_{T<\infty}] = \mathbb{E}[Z_T\mathbf{1}_{T<\infty}]$$
   *for any predictable stopping time $T$.*

**Remark 4.10.** Optional and predictable projections share some common properties with the conditional expectation: if $X$ is a measurable process and $Y$ is an optional (resp. predictable) bounded process, then

$$^o(XY) = Y\,^oX, \quad (\text{resp. } ^p(XY) = Y\,^pX).$$

Now, we state and prove an important theorem about the difference between the optional and predictable sigma fields.

**Theorem 4.11** ([69], p.174). *Let $\mathcal{I}$ be the sigma field generated by the processes $(M_t - M_{t-})$ where $M$ ranges through the bounded $(\mathcal{F}_t)$ martingales; then*

$$\mathcal{O} = \mathcal{P} \vee \mathcal{I}.$$

*In particular, if all $(\mathcal{F}_t)$ martingales are continuous, then $\mathcal{O} = \mathcal{P}$ and every stopping time is predictable.*

*Proof.* Since every optional projection is its own projection, it is enough to show that every optional projection is measurable with respect to $\mathcal{P} \vee \mathcal{I}$. But this is obvious for the processes $\mathbf{1}_{[0,u[}(t)H$ considered in the proof of Theorem 4.1 and an application of the monotone class theorem completes the proof. $\quad\square$

**Example 4.12.** In a Brownian filtration, every stopping time is predictable.

To conclude this paragraph, let us give an example from filtering theory.

**Example 4.13** ([69], exercise 5.15, p.175). Let $(\mathcal{F}_t)$ be a filtration, $B$ an $(\mathcal{F}_t)$ Brownian Motion and $h$ a bounded optional process. Let

$$Y_t = \int_0^t h_s ds + B_t;$$

if $\widehat{h}$ is the optional projection of $h$ with respect to $\left(\mathcal{F}_t^Y \equiv \sigma\left(Y_s,\ s \leq t\right)\right)$, then the process

$$N_t = Y_t - \int_0^t \widehat{h}_s ds$$

is an $\left(\mathcal{F}_t^Y\right)$ Brownian Motion. The process $N$ is called the innovation process.

### 4.2. Increasing processes and projections[3]

**Definition 4.14.** We shall call an increasing process a process $(A_t)$ which is nonnegative, $(\mathcal{F}_t)$ adapted, and whose paths are increasing and càdlàg.

**Remark 4.15.** A stochastic process which is nonnegative, and whose paths are increasing and càdlàg, but which is not $(\mathcal{F}_t)$ adapted, is called a *raw increasing process*.

Increasing processes play a central role in the general theory of stochastic processes: the main idea is to think of an increasing process as a random measure on $\mathbb{R}_+$, $dA_t(\omega)$, whose distribution function is $A_\bullet(\omega)$. That is why we shall make the convention that $A_{0-} = 0$, so that $A_0$ is the measure of $\{0\}$. The increasing process $A$ is called integrable if $\mathbb{E}[A_\infty] < \infty$.

A process which can be written as the difference of two increasing processes (resp. integrable increasing processes) is called a process of finite variation (resp. a process of integrable variation).

---

[3]We use again the convention $\mathcal{F}_{0-} = \mathcal{F}_0$.

Now, let $(A_t)$ be an increasing process and define its right continuous inverse:

$$C_t = \inf \left\{ s : \ A_s > t \right\}.$$

We have $C_{t-} = \inf \{s : \ A_s \geq t\}$. $C_t$ is the début of an optional set and hence it is a stopping time. Similarly, $C_{t-}$ is a stopping time which is predictable if $(A_t)$ is predictable (as the début of a right closed predictable set). Now, we recall the following time change formula ([27], p.132): for every nonnegative process $Z$, we have:

$$\int_0^\infty Z_s dA_s = \int_0^\infty Z_{C_{s-}} \mathbf{1}_{C_{s-} < \infty} ds \tag{4.1}$$

$$= \int_0^\infty Z_{C_s} \mathbf{1}_{C_s < \infty} ds. \tag{4.2}$$

With the time change formula, we can now prove the following theorem:

**Theorem 4.16.** *Let $X$ be a nonnegative measurable process and let $(A_t)$ be an increasing process. Then we have:*

$$\mathbb{E}\left[ \int_{[0,\infty[} X_s dA_s \right] = \mathbb{E}\left[ \int_{[0,\infty[} {}^o X_s dA_s \right],$$

*and if $(A_t)$ is predictable*

$$\mathbb{E}\left[ \int_{[0,\infty[} X_s dA_s \right] = \mathbb{E}\left[ \int_{[0,\infty[} {}^p X_s dA_s \right].$$

*Proof.* We shall prove the second equality which is more difficult to prove. From (4.1), we have:

$$\mathbb{E}\left[ \int_{[0,\infty[} X_s dA_s \right] = \mathbb{E}\left[ \int_0^\infty X_{C_{s-}} \mathbf{1}_{C_{s-} < \infty} ds \right] = \int_0^\infty ds \mathbb{E}\left[ X_{C_{s-}} \mathbf{1}_{C_{s-} < \infty} \right].$$

Similarly, we can show that:

$$\mathbb{E}\left[ \int_{[0,\infty[} {}^p X_s dA_s \right] = \int_0^\infty ds \mathbb{E}\left[ {}^p X_{C_{s-}} \mathbf{1}_{C_{s-} < \infty} \right].$$

Now, since $C_{s-}$ is a predictable stopping time, we have from the definition of predictable projections:

$$\mathbb{E}\left[ X_{C_{s-}} \mathbf{1}_{C_{s-} < \infty} \right] = \mathbb{E}\left[ {}^p X_{C_{s-}} \mathbf{1}_{C_{s-} < \infty} \right],$$

and this completes the proof of the theorem. $\qquad \square$

**Example 4.17.** One often uses the above theorem in the following way. Take $X_t(\omega) \equiv a(t) H(\omega)$, where $a(t)$ is a nonnegative Borel function and $H$ a nonnegative random variable. Let $(H_t)$ be a càdlàg version of the martingale $\mathbb{E}[H|\mathcal{F}_t]$. An application of Theorem 4.16 yields:

$$\mathbb{E}\left[H \int_{[0,\infty[} a(s) dA_s\right] = \mathbb{E}\left[\int_{[0,\infty[} H_s a(s) dA_s\right],$$

if $A$ is optional and

$$\mathbb{E}\left[H \int_{[0,\infty[} a(s) dA_s\right] = \mathbb{E}\left[\int_{[0,\infty[} H_{s-} a(s) dA_s\right],$$

if $A$ is predictable. In particular, for all uniformly $(\mathcal{F}_t)$ martingales $(M_t)$, we have (we assume $A_0 = 0$):

$$\mathbb{E}[M_\infty A_\infty] = \mathbb{E}\left[\int_{[0,\infty[} M_s dA_s\right],$$

if $A$ is optional and

$$\mathbb{E}[M_\infty A_\infty] = \mathbb{E}\left[\int_{[0,\infty[} M_{s-} dA_s\right],$$

if $A$ is predictable.

Now, we give a sort of converse to Theorem 4.16:

**Theorem 4.18** ([27], p.136). *Let $(A_t)$ be a raw increasing process which is further assumed to be integrable ($\mathbb{E}[A_\infty] < \infty$).*

1.  *If for all bounded and measurable processes $X$ we have*

$$\mathbb{E}\left[\int_{[0,\infty[} X_s dA_s\right] = \mathbb{E}\left[\int_{[0,\infty[} {}^o X_s dA_s\right],$$

   *then $A$ is optional.*

2.  *If for all bounded and measurable processes $X$ we have*

$$\mathbb{E}\left[\int_{[0,\infty[} X_s dA_s\right] = \mathbb{E}\left[\int_{[0,\infty[} {}^p X_s dA_s\right],$$

   *then $A$ is predictable.*

### *4.3. Random measures on $(\mathbb{R}_+ \times \Omega)$ and the dual projections*

**Definition 4.19.** We call $\mathbb{P}$-*measure* a bounded measure on the sigma field $\mathcal{B}(\mathbb{R}_+) \otimes \mathcal{F}$ (resp. $\mathcal{O}, \mathcal{P}$) which does not charge the $\mathbb{P}$-evanescent sets of $\mathcal{B}(\mathbb{R}_+) \otimes \mathcal{F}$ (resp. $\mathcal{O}, \mathcal{P}$).

We can construct a $\mathbb{P}$-measure $\mu$ in the following way: let $(A_t)$ be a raw increasing process which is integrable; for any bounded measurable process

$$\mu(X) \equiv \mathbb{E}\left[\int_{[0,\infty[} X_s dA_s\right].$$

In fact, quite remarkably, all $\mathbb{P}$-measures are of this form:

**Theorem 4.20** ([27], p.141). *Let $\mu$ be a nonnegative $\mathbb{P}$-measure (resp. $\mathbb{P}$-measure) on $\mathcal{B}(\mathbb{R}_+) \otimes \mathcal{F}$. Then, there exists a unique raw and integrable increasing process $(A_t)$ (resp. a raw process of integrable variation), up to indistinguishability, such that for any bounded process $X$*

$$\mu(X) = \mathbb{E}\left[\int_{[0,\infty[} X_s dA_s\right].$$

*We say that $A$ is the integrable increasing process (resp. the process of integrable variation) associated with $\mu$.*

*Furthermore, $(A_t)$ is optional (resp. predictable) if and only if for any bounded measurable process $X$:*

$$\mu(X) = \mu({}^o X); \quad \text{resp. } \mu(X) = \mu({}^p X).$$

**Definition 4.21.** A $\mathbb{P}$-measure is called optional (resp. predictable) if for any bounded measurable process $X$:

$$\mu(X) = \mu({}^o X); \quad \text{resp. } \mu(X) = \mu({}^p X).$$

**Remark 4.22.** From Theorem 4.18, the increasing process associated with an optional (resp. predictable) measure $\mu$ is optional (resp. predictable).

Now, we give some interesting consequences of Theorem 4.20. The first one is useful to prove the uniqueness of the Doob-Mayer decomposition of supermartingales.

**Theorem 4.23.** *Let $(A_t)$ and $(B_t)$ be two processes of integrable variation. If for any stopping time $T$*

$$\mathbb{E}\left[A_\infty - A_{T-}|\mathcal{F}_T\right] = \mathbb{E}\left[B_\infty - B_{T-}|\mathcal{F}_T\right] \text{ a.s. } (A_{0-} = B_{0-} = 0)$$

*then $A$ and $B$ are indistinguishable.*

*Similarly, if $(A_t)$ and $(B_t)$ are two predictable processes of integrable variation and if for any $t$*

$$\mathbb{E}\left[A_\infty - A_t|\mathcal{F}_t\right] = \mathbb{E}\left[B_\infty - B_t|\mathcal{F}_t\right] \text{ a.s.}$$

*then $A$ and $B$ are indistinguishable.*

*Proof.* We prove only the optional case, the proof for the predictable case being the same.

Let $\mu$ be the measure associated with $A - B$. The condition of the theorem entails that for any stopping time $T$:

$$\mu\left([T, \infty[\right) = 0,$$

and consequently, $\mu\left(U\right) = 0$, whenever $U$ is any finite and disjoint union of stochastic intervals of the form $[S, T[$. Hence $\mu$ is null on a Boole algebra which generates the optional sigma field and thus

$$\mu\left(^{o}X\right) = 0$$

for all bounded measurable process $X$, and from Theorem 4.20, $\mu$ is null on $\mathcal{B}\left(\mathbb{R}_{+}\right) \otimes \mathcal{F}$, and $A - B = 0$. $\quad\square$

**Remark 4.24.** The above theorem states that the processes $(A_{\infty} - A_{t-})$ and $(B_{\infty} - B_{t-})$ have the same optional projection. Moreover, as already observed for the projection theorems, one can replace the conditions of the theorem, in the optional case, with the unconditional form:

$$\mathbb{E}\left[A_{\infty} - A_{T-}\right] = \mathbb{E}\left[B_{\infty} - B_{T-}\right] \quad \text{for any stopping time } T.$$

Indeed, it suffices to consider the stopping times $T_H$, with $H \in \mathcal{F}_T$.

**Remark 4.25.** For the predictable case, one must be more cautious. First, we note that the condition of the theorem combined with the right continuity of the paths entail:

$$\mathbb{E}\left[A_{\infty} - A_T | \mathcal{F}_T\right] = \mathbb{E}\left[B_{\infty} - B_T | \mathcal{F}_T\right] \text{ a.s. for any stopping time } T.$$

From this we deduce now the unconditional form

$$\mathbb{E}\left[A_{\infty} - A_T\right] = \mathbb{E}\left[B_{\infty} - B_T\right] \quad \text{for any stopping time } T.$$

In this case, $(A_{\infty} - A_t)$ and $(B_{\infty} - B_t)$ have the same optional projection.

Now, we define projections of $\mathbb{P}$-measures:

**Definition 4.26.** Let $\mu$ be a $\mathbb{P}$-measure on $\mathcal{B}\left(\mathbb{R}_{+}\right) \otimes \mathcal{F}$. We call optional (resp. predictable) projection of $\mu$ the $\mathbb{P}$-measure $\mu^{o}$ (resp. $\mu^{p}$) defined on $\mathcal{B}\left(\mathbb{R}_{+}\right) \otimes \mathcal{F}$, by:

$$\mu^{o}\left(X\right) = \mu\left(^{o}X\right) = \quad (\text{resp. } \mu^{p}\left(X\right) = \mu\left(^{p}X\right)),$$

for all measurable and bounded processes $X$.

**Example 4.27.** Let $\mu$ be a nonnegative measure with a density $f$:

$$\mu\left(X\right) = \mathbb{E}\left[\int_{0}^{\infty} X_s\left(\omega\right) f\left(s, \omega\right) ds\right] \quad (f \geq 0, \ \mathbb{E}\left[\int_{0}^{\infty} f\left(s, \omega\right) ds\right] < \infty).$$

Then we have:

$$\mu^o\left(X\right) = \mu^p\left(X\right) = \mathbb{E}\left[\int_0^\infty X_s\left(\omega\right) g\left(s,\omega\right) ds\right],$$

where $g$ is any nonnegative measurable process such that

$$g\left(s,\bullet\right) = \mathbb{E}\left[f\left(s,\bullet\right)|\mathcal{F}_s\right]$$

for almost all $s$. For example $g$ can be the optional or predictable projection of $f$.

But usually, unlike the above example, $\mu^o$ (or $\mu^p$) is not associated with ${}^oA$ (or ${}^pA$); this leads us to the following fundamental definition of dual projections:

**Definition 4.28.** Let $(A_t)$ be an integrable raw increasing process. We call *dual optional projection* of $A$ the (optional) increasing process $(A_t^o)$ defined by:

$$\mathbb{E}\left[\int_{[0,\infty[} X_s dA_s^o\right] = \mathbb{E}\left[\int_{[0,\infty[} {}^oX_s dA_s\right],$$

for any bounded measurable $X$. We call *dual predictable projection* of $A$ the predictable increasing process $(A_t^p)$ defined by:

$$\mathbb{E}\left[\int_{[0,\infty[} X_s dA_s^p\right] = \mathbb{E}\left[\int_{[0,\infty[} {}^pX_s dA_s\right],$$

for any bounded measurable $X$.

**Remark 4.29.** The above definition extends in a straightforward way to processes of integrable variation.

**Remark 4.30.** Formally, the projection operation consists in defining conveniently the process $\mathbb{E}\left[A_t|\mathcal{F}_t\right]$, whereas the dual projection operation consists of defining the symbolic integral $\int_0^t \mathbb{E}\left[dA_s|\mathcal{F}_s\right]$. In particular, the dual projection of a bounded process needs not be bounded (for example, the dual projection of an honest time, as will be explained in a subsequent section).

Now, we shall try to compute the jumps of the dual projections.

**Proposition 4.31.** *Let $(A_t)$ and $(B_t)$ be two raw processes of integrable variation.*

1. *A and B have the same dual optional projection if and only if for any stopping time $T$:*

   $$\mathbb{E}\left[A_\infty - A_{T-}|\mathcal{F}_T\right] = \mathbb{E}\left[B_\infty - B_{T-}|\mathcal{F}_T\right] \text{ a.s. } (A_{0-} = B_{0-} = 0).$$

2. *A and B have the same dual predictable projection if and only if for every $t \geq 0$:*

   $$\mathbb{E}\left[A_\infty - A_t|\mathcal{F}_t\right] = \mathbb{E}\left[B_\infty - B_t|\mathcal{F}_t\right] \text{ a.s..}$$

*Proof.* It is exactly the same as the proof of Theorem 4.23. ∎

**Theorem 4.32.** *Let $(A_t)$ be a raw process of integrable variation.*

1. *Let $T$ be a stopping time. Then the jump of $A^o$ at $T$, $\Delta A_T^o$, is given by:*

$$\Delta A_T^o = \mathbb{E}\left[\Delta A_T | \mathcal{F}_T\right] \text{ a.s.,}$$

   *with the convention that $\Delta A_\infty^o = 0$.*
2. *Let $T$ be a predictable stopping time. Then the jump of $A^p$ at $T$, $\Delta A_T^p$, is given by:*

$$\Delta A_T^p = \mathbb{E}\left[\Delta A_T | \mathcal{F}_{T-}\right] \text{ a.s.,}$$

   *with the convention that $\Delta A_\infty^p = 0$.*

*Proof.* We only deal with the optional case (the predictable case can be dealt with similarly). From Proposition 4.31, we have:

$$\mathbb{E}\left[A_\infty - A_{T-}|\mathcal{F}_T\right] = \mathbb{E}\left[A_\infty^o - A_{T-}^o|\mathcal{F}_T\right].$$

Similarly, we have:

$$\mathbb{E}\left[A_\infty - A_T|\mathcal{F}_T\right] = \mathbb{E}\left[A_\infty^o - A_T^o|\mathcal{F}_T\right],$$

and consequently

$$\mathbb{E}\left[\Delta A_T|\mathcal{F}_T\right] = \mathbb{E}\left[\Delta A_T^o|\mathcal{F}_T\right].$$

Now, the result follows from the fact that $\Delta A_T^o$ is $\mathcal{F}_T$ measurable. ∎

Now we define predictable compensators.

**Definition 4.33.** Let $(A_t)$ be an optional process of integrable variation. The dual predictable projection of $A$, which we shall denote by $\widetilde{A}$, is also called the predictable compensator of $A$.

Why is $\widetilde{A}$ called compensator? From Proposition 4.31, we have:

$$A_t - \widetilde{A}_t = \mathbb{E}\left[A_\infty - \widetilde{A}_\infty|\mathcal{F}_t\right], \text{ a.s. and} \quad A_0 = \widetilde{A}_0.$$

Consequently, $A - \widetilde{A}$ is a martingale and $\widetilde{A}$ is the process that one has to subtract to $A$ to obtain a martingale which vanishes at 0. This is for example what one does to go from the Poisson process $N_t$ to the compensated Poisson process $N_t - \lambda t$. Another classical application is concerned with totally inaccessible stopping times:

**Proposition 4.34.** *Let $T$ be a stopping time which is strictly positive. Then the following are equivalent:*

1. *$T$ is totally inaccessible;*
2. *there exists a uniformly integrable martingale $(M_t)$, with $M_0 = 0$, which is continuous outside $[T]$, and such that $\Delta M_T = 1$ on $\{T < \infty\}$.*

*Proof.* $(2) \Rightarrow (1)$. We have to prove that $\mathbb{P}\left[T = S < \infty\right] = 0$ for any predictable stopping time $S$. We have:

$$\Delta M_S = \Delta M_S \mathbf{1}_{S < \infty} = \Delta M_T \mathbf{1}_{S=T<\infty} = \mathbf{1}_{S=T<\infty}.$$

Now, since $S$ is predictable,

$$\mathbb{E}\left[M_{S-}\right] = \mathbb{E}\left[M_0\right] = \mathbb{E}\left[M_S\right],$$

and hence

$$\mathbb{P}\left[T = S < \infty\right] = \mathbb{E}\left[\Delta M_S\right] = 0.$$

$(1) \Rightarrow (2)$. Let $A_t = \mathbf{1}_{T \leq t}$ and let $\widetilde{A}$ denote its predictable compensator. From Theorem 4.32, for all predictable stopping times $S$, we have:

$$\Delta \widetilde{A}_S = \mathbb{E}\left[\Delta A_S | \mathcal{F}_{S-}\right] = 0 \text{ a.s.}$$

since $\mathbb{P}\left[T = S < \infty\right] = 0$ ($T$ is predictable). Consequently, $\widetilde{A}$ is continuous. Now, it is easy to check that $M \equiv A - \widetilde{A}$ satisfies the required properties. $\qquad \square$

Now, we shall give the law of $\widetilde{A}_T = \widetilde{A}_\infty$.

**Proposition 4.35** (Azéma [2])**.** *Let $T$ be a finite, totally inaccessible stopping time. Let $\left(\widetilde{A}_t\right)$ be the predictable compensator of $\mathbf{1}_{T \leq t}$. Then $\widetilde{A}_T$ is exponentially distributed with parameter* $1$.

*Proof.* Let

$$M_t = \mathbf{1}_{T \leq t} - \widetilde{A}_t$$

be the martingale introduced in the proof of Proposition 4.34. We associate with $f$, a Borel bounded function, the stochastic integral $M_t^f \equiv \int_0^t f\left(A_s\right) dM_s$. We can compute $M_t^f$ more explicitly:

$$M_t^f = f\left(A_T\right) \mathbf{1}_{T \leq t} - \int_0^t f\left(A_s\right) dA_s = f\left(A_T\right) \mathbf{1}_{T \leq t} - F\left(A_t\right),$$

where $F\left(t\right) = \int_0^t ds f\left(s\right)$. An application of the optional stopping theorem yields $\mathbb{E}\left[M_\infty\right] = 0$, which implies:

$$\mathbb{E}\left[f\left(A_T\right)\right] = \mathbb{E}\left[F\left(A_T\right)\right].$$

As this last equality holds for every Borel bounded function, the law of $A_T$ must be the standard exponential law. $\qquad \square$

**Remark 4.36.** This above proposition has some nice applications to honest times that avoids stopping times in the theory of progressive enlargements of filtrations as we shall see later (see [42], [64]).

Now, we give two nice applications of dual projections. First, we refine Theorem 4.18:

**Theorem 4.37.** *Let* $(A_t)$ *be a raw increasing process such that* $\mathbb{E}[A_t] < \infty$ *for all* $t \geq 0$.

    *1. If for all bounded* $(\mathcal{F}_t)$ *martingales* $X$ *we have*

$$\mathbb{E}\left[\int_0^t X_s dA_s\right] = \mathbb{E}\left[X_\infty A_t\right], \text{ for all } t \geq 0$$

    *then* $A$ *is optional.*

    *2. If for all bounded* $(\mathcal{F}_t)$ *martingales* $X$ *we have*

$$\mathbb{E}\left[\int_0^t X_{s-} dA_s\right] = \mathbb{E}\left[X_\infty A_t\right], \text{ for all } t \geq 0$$

    *then* $A$ *is predictable.*

*Proof.* If $(X_t)$ is an $(\mathcal{F}_t)$ bounded martingale, then:

$$\mathbb{E}\left[\int_0^t X_s dA_s\right] = \mathbb{E}\left[\int_0^t X_s dA_s^o\right] = \mathbb{E}\left[X_t A_t^o\right] = \mathbb{E}\left[X_\infty A_t^o\right],$$

and similarly:

$$\mathbb{E}\left[\int_0^t X_{s-} dA_s\right] = \mathbb{E}\left[\int_0^t X_{s-} dA_s^p\right] = \mathbb{E}\left[X_t A_t^p\right] = \mathbb{E}\left[X_\infty A_t^p\right].$$

Consequently, under the hypothesis (1), one obtains:

$$\mathbb{E}\left[X A_t\right] = \mathbb{E}\left[X A_t^o\right], \quad \text{for every } X \in L^1,$$

which implies

$$A_t = A_t^o.$$

Likewise, under hypothesis (2), $A_t = A_t^p$. $\qquad\qquad\square$

    As a by product of Theorem 4.37, we can show the characterization of stopping times by Knight and Maisonneuve. We introduce the sigma field $(\mathcal{F}_\rho)$ associated with an arbitrary random time $\rho$, i.e. a nonnegative random variable:

$$\mathcal{F}_\rho = \sigma\left\{z_\rho, \ (z_t) \text{ any } (\mathcal{F}_t) \text{ optional process}\right\}.$$

**Theorem 4.38** (Knight-Maisonneuve [48]). *If for all uniformly integrable* $(\mathcal{F}_t)$ *martingales* $(M_t)$, *one has*

$$\mathbb{E}\left[M_\infty \mid \mathcal{F}_\rho\right] = M_\rho, \qquad \text{on } \{\rho < \infty\},$$

*then* $\rho$ *is a* $(\mathcal{F}_t)$ *stopping time (the converse is Doob's optional stopping theorem).*

*Proof.* For $t \geq 0$ we have

$$\mathbb{E}\left[M_\infty \mathbf{1}_{(\rho \leq t)}\right] = \mathbb{E}\left[M_\rho \mathbf{1}_{(\rho \leq t)}\right] = \mathbb{E}\left[\int_0^t M_s dA_s^\rho\right] = \mathbb{E}\left[M_\infty A_t^\rho\right].$$

Comparing the two extreme terms, we get

$$\mathbf{1}_{(\rho \leq t)} = A_t^\rho,$$

i.e $\rho$ is a $(\mathcal{F}_t)$ stopping time. $\qquad\square$

**Remark 4.39.** The result of Knight and Maisonneuve suggests the following questions:

- How may $\mathbb{E}\left[M_\infty \mid \mathcal{F}_\rho\right]$ on the one hand, and $M_\rho$ on the other hand, differ for a non stopping time $\rho$? The reader can refer to [6], [81] or [64] for some answers.
- Given an arbitrary random time $\rho$, is it possible to characterize the set of martingales which satisfy (2.1)? (see [6], [81] or [64]).

Now we shall see how an application of dual projections and their simple properties gives a simple proof of a multidimensional extension of Lévy's arc sine law. The results that follow are borrowed from [62]. Consider $(R_t, L_t)$, where $(R_t)$ is a Bessel process of dimension $\delta \equiv 2(1 - \mu) \in (0, 2)$, starting from 0, and $(L_t)$ a normalization of its local time at level zero (see [62] for more precisions on this normalization).

**Lemma 4.40** ([62]). *Let*

$$g_\mu \equiv \sup \{t \leq 1 : \ R_t = 0\}.$$

*Then the dual predictable projection $A_t^{g_\mu}$ of $\mathbf{1}_{(g_\mu \leq t)}$ is:*

$$A_t^{g_\mu} = \frac{1}{2^\mu \Gamma(1 + \mu)} \int_0^{t \wedge 1} \frac{dL_u}{(1 - u)^\mu},$$

*i.e. for every nonnegative predictable process $(x_t)$,*

$$\mathbb{E}\left[x_{g_\mu}\right] = \frac{1}{2^\mu \Gamma(1 + \mu)} \mathbb{E}\left[\int_0^1 dL_u \frac{x_u}{(1 - u)^\mu}\right].$$

**Remark 4.41.** The random time $g_\mu$ is not a stopping time; it is a typical example of an honest time (we shall define these times in a subsequent section).

Now we give a result which was obtained by Barlow, Pitman and Yor, using excursion theory. The proof we give is borrowed from [62].

**Proposition 4.42** ([14],[62]). *The variable $g_\mu$ follows the law:*

$$\mathbb{P}(g_\mu \in dt) = \frac{\sin(\mu\pi)}{\pi} \frac{dt}{t^{1-\mu}(1-t)^\mu}, \ 0 < t < 1,$$

*i.e. the Beta law with parameters $(\mu, 1 - \mu)$. In particular, $\mathbb{P}(g \in dt) = \dfrac{1}{\pi} \dfrac{dt}{\sqrt{t(1-t)}}$,*

*i.e. $g$ is arc sine distributed ([49]).*

*Proof.* From Lemma 4.40, for every Borel function $f : [0, 1] \to \mathbb{R}_+$, we have:

$$\mathbb{E}\left[f\left(g_\mu\right)\right] = \frac{1}{2^\mu \mu \Gamma(\mu)} \mathbb{E}\left[\int_0^1 dL_u \frac{f(u)}{(1-u)^\mu}\right] = \frac{1}{2^\mu \mu \Gamma(\mu)} \int_0^1 du \, \mathbb{E}\left[L_u\right] \frac{f(u)}{(1-u)^\mu}. \tag{4.3}$$

By the scaling property of $(L_t)$,

$$\mathbb{E}\left[L_u\right] = u^\mu \mathbb{E}\left[L_1\right].$$

Moreover, by definition of $(L_t)$ (see [62] or [21]),

$$\mathbb{E}\left[L_1\right] = \mathbb{E}\left[R_1^{2\mu}\right];$$

since $R_1^2$ is distributed as $2\,gam(1-\mu)$, where $gam(1-\mu)$ denotes a random variable which follows the gamma law with parameter $\mu$, we have

$$\mathbb{E}\left[R_1^{2\mu}\right] = \frac{2^\mu}{\Gamma(1-\mu)}.$$

Now, plugging this in (4.3) yields:

$$\mathbb{E}\left[f\left(g_\mu\right)\right] = \frac{1}{\Gamma(\mu)\Gamma(1-\mu)} \int_0^1 du \frac{f(u)}{t^{1-\mu}(1-u)^\mu}.$$

To conclude, it suffices to use the duplication formula for the Gamma function ([1]):

$$\Gamma(\mu)\Gamma(1-\mu) = \frac{\pi}{\sin(\mu\pi)}.$$

<div style="text-align: right;">□</div>

We shall conclude this section giving a result which is useful in stochastic integration. Before, we need to introduce the notion of locally integrable increasing process.

**Definition 4.43.** A raw increasing process, such that $A_0 = 0$, is said to be locally integrable if there exists an increasing sequence of stopping times $(T_n)$, such that $\lim_{n\to\infty} T_n = \infty$, and

$$\mathbb{E}\left[A_{T_n}\right] < \infty.$$

**Theorem 4.44.** *[[27], Theorem 80, p.153]*

  1. *Every predictable process of finite variation is locally integrable.*
  2. *An optional process $A$ of finite variation is locally integrable if and only if there exists a predictable process $\widetilde{A}$ of finite variation such that $(A - \widetilde{A})$ is a local martingale which vanishes at 0. When it exists, $\widetilde{A}$ is unique. We say that $\widetilde{A}$ is the predictable compensator of $A$.*

**Remark 4.45.** We shall see an application of this result to the existence of the bracket $\langle M \rangle$ of a local martingale in next section.

## 5. The Doob-Meyer decomposition and multiplicative decompositions

This section is devoted to a fundamental result in Probability Theory: the Doob-Meyer decomposition theorem. This result is now standard and a proof of it can be found for example in [27, 68, 70].

**Theorem 5.1** (Doob-Meyer decomposition). *An adapted càdlàg process $Y$ is a submartingale of class $(D)$ null at $0$ if and only if $Y$ may be written:*

$$Y_t = M_t + A_t \tag{5.1}$$

*where $M$ is a uniformly integrable martingale null at $0$ and $A$ a predictable integrable increasing process null at $0$. Moreover, the decomposition (5.1) is unique.*

The following theorem gives a necessary and sufficient condition for $A$ to be continuous.

**Proposition 5.2.** *Let $Y = M + A$ be a submartingale of class $(D)$. The compensator $A$ of $Y$ is continuous if and only if for every predictable stopping time $T$:*

$$\mathbb{E}\left[Y_T\right] = \mathbb{E}\left[Y_{T-}\right].$$

*Proof.* Let $(T_n)$ be a sequence of stopping times that announce $T$. Since $Y$ is of class $(D)$, we have:

$$\lim_{n \to \infty} \mathbb{E}\left[Y_\infty - Y_{T_n}\right] = \lim_{n \to \infty} \mathbb{E}\left[A_\infty - A_{T_n}\right],$$

which is:

$$\mathbb{E}\left[Y_\infty - Y_{T-}\right] = \mathbb{E}\left[A_\infty - A_{T-}\right].$$

Now, we also have $\mathbb{E}\left[Y_\infty - Y_T\right] = \mathbb{E}\left[A_\infty - A_T\right]$. Consequently, we have:

$$\mathbb{E}\left[Y_T - Y_{T-}\right] = \mathbb{E}\left[A_T - A_{T-}\right],$$

and the result of the proposition follows easily. □

Now we give the local form of the Doob-Meyer decomposition; for this we need the following lemma:

**Lemma 5.3** ([70], p.374). *A submartingale $Y$ is locally of class $(D)$ (i.e. there exists a sequence of stopping times $T_n$ such that $T_n \to \infty$ as $n \to \infty$ and such that the stopped process $(Y_{t \wedge T_n})$ is a submartingale of class $(D)$ for every $n$).*

**Theorem 5.4.** *Let $Y$ be a local submartingale. Then $Y$ may be written uniquely in the form:*

$$Z = Z_0 + M + A$$

*where $M$ is a local martingale null at $0$ and $A$ is a predictable increasing process null at $0$.*

**Remark 5.5.** All the previous results were stated for submartingales but of course they also hold for supermartingales. Given a supermartingale $Z$, it suffices to consider the submartingale $Y = -Z$ to obtain the decomposition: $Z = M - A$.

To conclude the discussion on the Doob-Meyer decomposition, we give a result on the existence of the bracket $\langle X \rangle$ of a local martingale, which follows easily from the above theorem and Theorem 4.44.

**Theorem 5.6.** *Let $X$ be a local martingale null at $0$. The following are equivalent:*

1. *There exists a unique predictable increasing process $(\langle X \rangle_t)$ null at $0$ such that $X^2 - \langle X \rangle_t$ is a local martingale;*
2. *$X$ is locally $L^2$ bounded;*
3. *$[X]$ is locally integrable.*

*When one of these condition holds, we have:*

$$\langle X \rangle_t = \widetilde{[X]}_t.$$

## 6. Multiplicative decompositions

From this section on, we shall deal with more recent and less known aspects of the general theory of stochastic processes and stochastic calculus.

We start with a nice result about the existence of multiplicative decompositions for nonnegative submartingales and supermartingales. These decompositions have appeared much less useful than the additive decompositions. However, we mention some of them and give a nice application. The reader should refer to [36, 55, 3, 37, 65] for more details. We will also give a very elegant application of multiplicative decompositions to the theory of enlargements of filtrations in Section 8.

The first result in this direction is due to Itô and Watanabe and deals with nonnegative supermartingales. Let $(Z_t)$ be a nonnegative càdlàg supermartingale, and define

$$T_0 = \inf \{t :\ Z_t = 0\}.$$

It is well known that $Z$ is null on $[T, \infty[$.

**Theorem 6.1** (Itô-Watanabe [36]). *Let $(Z_t)$ be a nonnegative càdlàg supermartingale such that $\mathbb{P}(T_0 > 0) = 1$. Then $Z$ can be factorized as:*

$$Z_t = Z_t^{(0)} Z_t^{(1)},$$

*with a positive local martingale $Z_t^{(0)}$ and a decreasing process $Z_t^{(1)}$ ($Z_0^{(1)} = 1$). If there are two such factorizations, then they are identical in $[0, T_0[$.*

**Remark 6.2.** We shall see in subsequent sections how a refinement of this decomposition in the special case of the Azéma's supermartingales associated with honest times leads to some nice results on enlargements of filtration ([63]).

Now, we introduce a remarkable family of continuous local submartingales, which we shall call of class $(\Sigma_c)$ (the subscript c stands for continuous) and which appear under different forms in Probability Theory (in the study of local times, of supremum of local martingales, in the balayage,etc.). The reader can refer to [69, 61, 60, 65] for more references and details[4]

**Definition 6.3.** Let $(X_t)$ be a positive local submartingale, which decomposes as:

$$X_t = N_t + A_t. \tag{6.1}$$

We say that $(X_t)$ is of class $(\Sigma_c)$ if:

1. $(N_t)$ is a continuous local martingale, with $N_0 = 0$;
2. $(A_t)$ is a continuous increasing process, with $A_0 = 0$;
3. the measure $(dA_t)$ is carried by the set $\{t : X_t = 0\}$.

If additionally, $(X_t)$ is of class $(D)$, we shall say that $(X_t)$ is of class $(\Sigma_c D)$. Now, if in (1) $N$ is only assumed to be càdlàg, we say that $(X_t)$ is of class $(\Sigma)$, dropping the subscript $c$. Similarly we define the class $(\Sigma D)$.

**Example 6.4.** The absolute value of a local martingale $|M_t|$, $S_t - M_t$, where $S_t \equiv \sup_{u \leq t} M_u$, $\alpha M_t^+ + \beta M_t^-$ with $\alpha > 0$, $\beta > 0$, or $R_t^{2\mu}$ where $(R_t)$ is a Bessel process, starting from 0, of dimension $2(1 - \mu)$, with $\mu \in (0, 1)$, are typical examples of processes of class $(\Sigma_c)$.

We shall use a decomposition result for nonnegative local submartingales to obtain a characterization of submartingales of class $(\Sigma_c)$ among other local submartingales. Let $Y$ be a continuous nonnegative local submartingale. It is in general impossible to get a multiplicative decomposition of the form $Y_t = Y_t^{(0)} Y_t^{(1)}$, with a positive local martingale $Y_t^{(0)}$ and an increasing process $Y_t^{(1)}$ ($Y_0^{(1)} = 1$). Indeed, it is well known that once a nonnegative local martingale is equal to zero, then it remains null, while this is not true for nonnegative submartingales (consider for example $|B_t|$, the absolute value of a standard Brownian Motion). Hence we shall use the following form of multiplicative decomposition:

**Proposition 6.5** ([65])**.** *Let $(Y_t)_{t \geq 0}$ be a continuous nonnegative local submartingale such that $Y_0 = 0$. Consider its Doob-Meyer decomposition:*

$$Y_t = m_t + \ell_t. \tag{6.2}$$

*The local submartingale $(Y_t)_{t \geq 0}$ then admits the following multiplicative decomposition:*

$$Y_t = M_t C_t - 1, \tag{6.3}$$

*where $(M_t)_{t \geq 0}$ is a continuous local martingale, which is strictly positive, with $M_0 = 1$ and where $(C_t)_{t \geq 0}$ is an increasing continuous and adapted process,*

---

[4]In [60, 65], this class is simply called $(\Sigma)$ while in [61] the class $(\Sigma)$ is bigger in the sense that $N$ can be càdlàg.

*with $C_0 = 1$. The decomposition is unique and the processes $C$ and $M$ are given by the explicit formulae:*

$$C_t = \exp\left(\int_0^t \frac{d\ell_s}{1 + Y_s}\right), \tag{6.4}$$

*and*

$$M_t = (Y_t + 1)\exp\left(-\int_0^t \frac{d\ell_s}{1 + Y_s}\right) \tag{6.5}$$

$$= \exp\left(\int_0^t \frac{dm_s}{1 + Y_s} - \frac{1}{2}\int_0^t \frac{d\langle m \rangle_s}{(1 + Y_s)^2}\right). \tag{6.6}$$

**Remark 6.6.** It is possible to find necessary and sufficient conditions on $M$ and $C$ for $Y$ to be of class $(D)$ (see [65]).

Now, if we want $Y$ to be of class $(\Sigma_c)$, Proposition 6.5 takes the following more precise form:

**Proposition 6.7** ([65])**.** *Let $(X_t = N_t + A_t, \ t \geq 0)$ be a nonnegative, continuous local submartingale with $X_0 = 0$. Then, the following are equivalent:*

1. *$(X_t, \ t \geq 0)$ is of class $(\Sigma_c)$, i.e. $(dA_t)$ is carried by the set $\{t : \ X_t = 0\}$.*
2. *There exists a strictly positive, continuous local martingale $(M_t)$, with $M_0 = 1$, such that:*

$$X_t = \frac{M_t}{I_t} - 1, \tag{6.7}$$

   *where*

$$I_t = \inf_{s \leq t} M_s.$$

*The local martingale $(M_t)$ is given by:*

$$M_t = (1 + X_t)\exp\left(-A_t\right). \tag{6.8}$$

Now, we can give a nice characterization of the local submartingales of the class $(\Sigma_c)$ in terms of frequency of vanishing. More precisely, let $(M_t)$ be a strictly positive and continuous local martingale with $M_0 = 1$ and denote by $\mathcal{E}(M)$ the set of all nonnegative local submartingales with the same martingale part $M$ in their multiplicative decomposition (6.3). Then the following holds:

**Corollary 6.8** ([65])**.** *Let*

$$Y^\star = \frac{M_t}{I_t} - 1.$$

*Then, $(Y_t^\star)$ is in $\mathcal{E}(M)$ and it is the smallest element of $\mathcal{E}(M)$ in the sense that:*

$$\forall Y \in \mathcal{E}(M), \ Y^\star \leq Y.$$

*Consequently, $(Y_t^\star)$ has more zeros than any other local submartingale of $\mathcal{E}(M)$.*

*Proof.* It suffices to note that any element $Y \in \mathcal{E}(M)$ decomposes as $Y_t = M_t C_t - 1$. Since $Y$ must be nonnegative, we must have:

$$M_t \geq \frac{1}{C_t}.$$

But $\frac{1}{C_t}$ is decreasing, hence we have:

$$\frac{1}{C_t} \leq I_t,$$

and this proves the Corollary. $\qquad\qquad\qquad\qquad\qquad\qquad\qquad\qquad\square$

## 7. Some hidden martingales

In this section, we illustrate the power of martingale methods by extending some (well) known results for the Brownian Motion to larger classes of stochastic processes which do not in general enjoy scaling or Markov properties.

First, we need the following lemma which we will not prove since its proof is similar to the proof of Proposition 3.6 (it is in fact much simpler, see [69], Proposition 3.5, p.70).

**Lemma 7.1.** *Let $(X_t)$ be a càdlàg, adapted and bounded process. Then $(X_t)$ is a martingale if and only if for any bounded stopping time $T$:*

$$\mathbb{E}[X_T] = \mathbb{E}[X_0].$$

Now, we state and prove with stochastic calculus arguments a characterization of predictable increasing processes among optional processes.

**Theorem 7.2.** *Let $(A_t)$ be an increasing optional process such that:*

$$\mathbb{E}[A_\infty] < \infty.$$

*Then $(A_t)$ is predictable if and only if for every bounded martingale $(X_t)$,*

$$\mathbb{E}\left[\int_0^\infty X_s dA_s\right] = \mathbb{E}\left[\int_0^\infty X_{s-} dA_s\right]. \tag{7.1}$$

*Proof.* We first note that (7.1) is equivalent to:

$$\mathbb{E}\left[\sum_{s \geq 0} (\Delta X_s)(\Delta A_s)\right] = 0,$$

which by stopping is equal to:

$$\mathbb{E}\left[\sum_{s \leq T} (\Delta X_s)(\Delta A_s)\right] = 0, \quad \text{for all bounded stopping times } T.$$

From Lemma 7.1, this is equivalent to:

$$\sum_{s \leq \bullet} (\Delta X_s)(\Delta A_s) \quad \text{is a martingale.} \tag{7.2}$$

Now, if $A$ is predictable, it is a classical result that ([76])

$$\sum_{s \leq \bullet} (\Delta X_s)(\Delta A_s) = \int_0^{\bullet} (\Delta A_s) \, dM_s,$$

and hence Theorem 7.2 follows.

Now, conversely, if (7.1) (or equivalently 7.2) is satisfied, then we can integrate any bounded predictable process $H$ with respect to this martingale, and the result is still a martingale (starting from 0). Thus we have obtained:

$$\mathbb{E}\left[\int_0^{\infty} H_s X_s \, dA_s\right] = \mathbb{E}\left[\int_0^{\infty} H_s X_{s-} \, dA_s\right] = \mathbb{E}\left[\int_0^{\infty} H_s X_{s-} \, dA_s^p\right].$$

We now take $X \equiv 1$ and we obtain, for any bounded predictable process $H$:

$$\mathbb{E}\left[\int_0^{\infty} H_s \, dA_s\right] = \mathbb{E}\left[\int_0^{\infty} H_s \, dA_s^p\right].$$

Hence $A = A^p$ and $A$ is predictable. $\square$

**Remark 7.3.** Since $(X_t)$ is the optional projection of $X_\infty$, then $\mathbb{E}\left[\int_0^{\infty} X_s \, dA_s\right] = \mathbb{E}[X_\infty A_\infty]$ and consequently (7.1) is equivalent to:

$$\mathbb{E}[X_\infty A_\infty] = \mathbb{E}\left[\int_0^{\infty} X_{s-} \, dA_s\right].$$

**Corollary 7.4.** *Let $T$ be a stopping time. $T$ is a predictable stopping time if and only if, for all bounded càdlàg martingales $(X_t)$*

$$\mathbb{E}[X_T \mathbf{1}_{T<\infty}] = \mathbb{E}[X_{T-} \mathbf{1}_{T<\infty}].$$

Now, we give an illustration of the power of martingale methods. We start with two similar results on Brownian Motion, obtained with excursion theory and Markov processes methods, and then we show how these results can be extended to a much wider class of processes, using martingale methods. This is also here the opportunity for us to review some classical martingale techniques. We start with a definition, which extends the class $(\Sigma_c)$ to a wider class $(\Sigma)$ which also contains some discontinuous martingales, such as Azéma's second martingale or its generalization (see [62]).

**Definition 7.5.** A nonnegative local submartingale $(X_t)$ of class $(\Sigma_c)$ is said to be of class $(\Sigma)$ if in condition (1) of definition 6.3 the local martingale $N$ can be chosen càdlàg (instead of continuous). Consequently, $(\Sigma_c) \subset (\Sigma)$.

**Remark 7.6.** In addition to the examples given after definition 6.3, we can give the following ones:

- Let $(M_t)$ be a local martingale (starting from 0) with only negative jumps and let $S_t \equiv \sup_{u \leq t} M_u$; then

$$X_t \equiv S_t - M_t$$

is of class $(\Sigma)$. In this case, $X$ has only positive jumps.

- Let $(R_t)$ be a Bessel process (starting from 0) of dimension $2(1 - \mu)$, with $\mu \in (0, 1)$. Define:

$$g_\mu(t) \equiv \sup \{u \leq t : R_u = 0\}.$$

In the filtration $\mathcal{G}_t \equiv \mathcal{F}_{g_\mu(t)}$ of the zeros of the Bessel process $R$, the stochastic process:

$$X_t \equiv (t - g_\mu(t))^\mu,$$

is a submartingale of class $(\Sigma)$ whose increasing process in its Doob-Meyer decomposition is given by:

$$A_t \equiv \frac{1}{2^\mu \Gamma(1 + \mu)} L_t(R),$$

where as usual $\Gamma$ stands for Euler's gamma function. Recall that $\mu \equiv \frac{1}{2}$ corresponds to the absolute value of the standard Brownian Motion; thus for $\mu \equiv \frac{1}{2}$ the above result leads to nothing but the celebrated second Azéma's martingale ($X_t - A_t$, see [11, 81]). In this example, $X$ has only negative jumps.

The local submartingales of the class $(\Sigma)$ have the following nice characterization based on stochastic calculus:

**Theorem 7.7** ([61])**.** *The following are equivalent:*

1. *The local submartingale $(X_t)$ is of class $(\Sigma)$;*
2. *There exists an increasing, adapted and continuous process $(C_t)$ such that for every locally bounded Borel function $f$, and $F(x) \equiv \int_0^x f(z)\, dz$, the process*

$$F(C_t) - f(C_t) X_t$$

*is a local martingale. Moreover, in this case, $(C_t)$ is equal to $(A_t)$, the increasing process of $X$.*

*Proof.* $(1) \Longrightarrow (2)$. First, let us assume that $f$ is $\mathcal{C}^1$ and let us take $C_t \equiv A_t$. An integration by parts yields:

$$
\begin{aligned}
f(A_t) X_t &= \int_0^t f(A_u)\, dX_u + \int_0^t f'(A_u)\, X_u dA_u \\
&= \int_0^t f(A_u)\, dN_u + \int_0^t f(A_u)\, dA_u + \int_0^t f'(A_u)\, X_u dA_u.
\end{aligned}
$$

Since $(dA_t)$ is carried by the set $\{t:\ X_t = 0\}$, we have $\int_0^t f'(A_u) X_u dA_u = 0$. As $\int_0^t f(A_u)\, dA_u = F(A_t)$, we have thus obtained that:

$$F(A_t) - f(A_t) X_t = -\int_0^t f(A_u)\, dN_u, \tag{7.3}$$

and consequently $(F(A_t) - f(A_t) X_t)$ is a local martingale. The general case when $f$ is only assumed to be locally bounded follows from a monotone class argument and the integral representation (7.3) is still valid.

$(2) \implies (1)$. First take $F(a) = a$; we then obtain that $C_t - X_t$ is a local martingale. Hence the increasing process of $X$ in its Doob-Meyer decomposition is $C$, and $C = A$. Next, we take: $F(a) = a^2$ and we get:

$$A_t^2 - 2A_t X_t$$

is a local martingale. But

$$A_t^2 - 2A_t X_t = 2\int_0^t A_s\,(dA_s - dX_s) - 2\int_0^t X_s dA_s$$

Hence, we must have:

$$\int_0^t X_s dA_s = 0.$$

Thus $dA_s$ is carried by the set of zeros of $X$. $\qquad\square$

Now, we state and prove the so called Doob's maximal identity, which is obtained as an easy application of Doob's optional stopping theorem, but which has many nice and deep applications (see [63]).

**Lemma 7.8** (Doob's maximal identity). *Let $(M_t)$ be a positive local martingale which satisfies:*

$$M_0 = x,\ x > 0;\ \lim_{t\to\infty} M_t = 0.$$

*If we note*

$$S_t \equiv \sup_{u\le t} M_u,$$

*and if $S$ is continuous, then for any $a > 0$, we have:*

*1.*

$$\mathbf{P}(S_\infty > a) = \left(\frac{x}{a}\right) \wedge 1. \tag{7.4}$$

*Hence, $\dfrac{x}{S_\infty}$ is a uniform random variable on $(0,1)$.*

*2. For any stopping time $T$:*

$$\mathbf{P}\left(S^T > a \mid \mathcal{F}_T\right) = \left(\frac{M_T}{a}\right) \wedge 1, \tag{7.5}$$

*where*

$$S^T = \sup_{u\ge T} M_u.$$

Hence $\dfrac{M_T}{S^T}$ is also a uniform random variable on $(0,1)$, independent of $\mathcal{F}_T$.

*Proof.* Formula (7.5) is a consequence of (7.4) when applied to the martingale $(M_{T+u})_{u \geq 0}$ and the filtration $(\mathcal{F}_{T+u})_{u \geq 0}$. Formula (7.4) itself is obvious when $a \leq x$, and for $a > x$, it is obtained by applying Doob's optional stopping theorem to the local martingale $(M_{t \wedge T_a})$, where $T_a = \inf \{u \geq 0 : M_u \geq a\}$. $\square$

Now let us mention a result of Knight ([46, 47]) which motivated our study:

**Proposition 7.9** (Knight [47])**.** *Let $(B_t)$ denote a standard Brownian Motion, and $S$ its supremum process. Then, for $\varphi$ a nonnegative Borel function, we have:*

$$\mathbb{P} \left( \forall t \geq 0, \ S_t - B_t \leq \varphi\left(S_t\right) \right) = \exp \left( - \int_0^\infty \frac{dx}{\varphi\left(x\right)} \right).$$

*Furthermore, if we let $T_x$ denote the stopping time:*

$$T_x = \inf \{t \geq 0 : \ S_t > x\} = \inf \{t \geq 0 : \ B_t > x\},$$

*then for any nonnegative Borel function $\varphi$, we have:*

$$\mathbb{P} \left( \forall t \leq T_x, \ S_t - B_t \leq \varphi\left(S_t\right) \right) = \exp \left( - \int_0^x \frac{dx}{\varphi\left(x\right)} \right).$$

Now, we give a more general result, which can also be applied to some discontinuous processes such as the one parameter generalizations of Azéma's second submartingale:

**Theorem 7.10** ([61])**.** *Let $X$ be a local submartingale of the class $(\Sigma)$, with only negative jumps, such that $\lim_{t \to \infty} A_t = \infty$. Define $(\tau_u)$ the right continuous inverse of $A$:*

$$\tau_u \equiv \inf \{t : \ A_t > u\}.$$

*Let $\varphi : \mathbb{R}_+ \to \mathbb{R}_+$ be a Borel function. Then, we have the following estimates:*

$$\mathbb{P} \left( \exists t \geq 0, \ X_t > \varphi\left(A_t\right) \right) = 1 - \exp \left( - \int_0^\infty \frac{dx}{\varphi\left(x\right)} \right), \qquad (7.6)$$

*and*

$$\mathbb{P} \left( \exists t \leq \tau_u, \ X_t > \varphi\left(A_t\right) \right) = 1 - \exp \left( - \int_0^u \frac{dx}{\varphi\left(x\right)} \right). \qquad (7.7)$$

*Proof.* The proof is based on Theorem 7.7 and Lemma 7.8. We shall first prove equation (7.6), and for this, we first note that we can always assume that $\dfrac{1}{\varphi}$ is bounded and integrable. Indeed, let us consider the event

$$\Delta_\varphi \equiv \{\exists t \geq 0, \ X_t > \varphi\left(A_t\right)\}.$$

Now, if $(\varphi_n)_{n \geq 1}$ is a decreasing sequence of functions with limit $\varphi$, then the events $(\Delta_{\varphi_n})$ are increasing, and $\bigcup_n \Delta_{\varphi_n} = \Delta_\varphi$. Hence, by approximating $\varphi$ from above, we can always assume that $\frac{1}{\varphi}$ is bounded and integrable.

Now, let

$$F(x) \equiv 1 - \exp\left(-\int_x^\infty \frac{dz}{\varphi(z)}\right);$$

its Lebesgue derivative $f$ is given by:

$$f(x) = \frac{-1}{\varphi(x)} \exp\left(-\int_x^\infty \frac{dz}{\varphi(z)}\right) = \frac{-1}{\varphi(x)}(1 - F(x)).$$

Now, from Theorem 7.7, $(M_t \equiv F(A_t) - f(A_t) X_t)$, which is also equal to $F(A_t) + X_t \frac{1}{\varphi(A_t)}(1 - F(A_t))$, is a positive local martingale (whose supremum is continuous since $(M_t)$ has only negative jumps), with $M_0 = 1 - \exp\left(-\int_0^\infty \frac{dx}{\varphi(x)}\right)$. Moreover, as $(M_t)$ is a positive local martingale, it converges almost surely as $t \to \infty$. Let us now consider $M_{\tau_u}$:

$$M_{\tau_u} = F(u) - f(u) X_{\tau_u}.$$

But since $(dA_t)$ is carried by the zeros of $X$ and since $\tau_u$ corresponds to an increase time of $A$, we have $X_{\tau_u} = 0$. Consequently,

$$\lim_{u \to \infty} M_{\tau_u} = \lim_{u \to \infty} F(u) = 0,$$

and hence

$$\lim_{u \to \infty} M_u = 0.$$

Now let us note that if for a given $t_0 < \infty$, we have $X_{t_0} > \varphi(A_{t_0})$, then we must have:

$$M_{t_0} > F(A_{t_0}) - f(A_{t_0}) \varphi(A_{t_0}) = 1,$$

and hence we easily deduce that:

$$
\begin{aligned}
\mathbb{P}\left(\exists t \geq 0, \ X_t > \varphi(A_t)\right) &= \mathbb{P}\left(\sup_{t \geq 0} M_t > 1\right) \\
&= \mathbb{P}\left(\sup_{t \geq 0} \frac{M_t}{M_0} > \frac{1}{M_0}\right) \\
&= M_0,
\end{aligned}
$$

where the last equality is obtained by an application of Doob's maximal identity (Lemma 7.8).

To obtain the second identity of the Theorem, it suffices to replace $\varphi$ by the function $\varphi_u$ defined as:

$$\varphi_u(x) = \begin{cases} \varphi(x) & \text{if } x < u \\ \infty & \text{otherwise.} \end{cases}$$

$\square$

**Remark 7.11.** The estimate of Knight is a consequence of Theorem 7.10, with $X_t = S_t - B_t$, and $A_t = S_t$. For applications of Theorem 7.10 to the processes $(t - g_\mu(t))^\mu$ introduced in Remark 7.6 and to the Skorokhod's stopping problem, see [61].

## 8. General random times, their associated $\sigma$-fields and Azéma's supermartingales

The role of stopping times in Probability Theory is fundamental and there are myriads of applications of the optional stopping theorems (2.1) and (2.2). However, it often happens that one needs to work with random times which are not stopping times: for example, in mathematical finance, in the modeling of default times (see [30] or [40] for an account and more references) or in insider trading models ([33]); in Markov Processes theory (see [28]); in the characterization of the set of zeros of continuous martingales ([12]), in path decomposition of some diffusions (see [73], [57, 58], [42] or [63]); in the study of Strong Brownian Filtrations (see [18]), etc. One of the aims of this essay is to go beyond the classical (yet important) concept of stopping times. One of the main tools to study random times which are not stopping times is the theory of progressive enlargements (or expansions) of filtrations. This section is devoted to important definitions and results from the general theory of stochastic processes which are useful to develop the theory of progressive enlargements of filtrations.

We first give the definition of BMO and $\mathcal{H}^1$ spaces which we shall use in the sequel. For more details, the reader can refer to [27].

Let $\left(\Omega, \mathcal{F}, (\mathcal{F}_t)_{t \geq 0}, \mathbb{P}\right)$ be a filtered probability space. We recall that the space $\mathcal{H}^1$ is the Banach space of (càdlàg) $(\mathcal{F}_t)$-martingales $(M_t)$ such that

$$\|M\|_{\mathcal{H}^1} = \mathbb{E}\left[\sup_{t \geq 0} |M_t|\right] < \infty.$$

The space of BMO martingales is the Banach space of (càdlàg) square integrable $(\mathcal{F}_t)$-martingales $(Y_t)$ which satisfy

$$\|Y\|_{BMO}^2 = \operatorname{esssup}_T \mathbb{E}\left[(Y_\infty - Y_{T-})^2 \mid \mathcal{F}_T\right] < \infty$$

where $T$ ranges over all $(\mathcal{F}_t)$-stopping times. It is a very nice result of Meyer that the dual of the space $\mathcal{H}^1$ is the space BMO ([52]).

### 8.1. Arbitrary random times and some associated sigma fields

**Definition 8.1.** A random time $\rho$ is a nonnegative random variable $\rho : (\Omega, \mathcal{F}) \to [0, \infty]$.

There are essentially two classes of random times that have been studied in detail which are not stopping times: ends of optional or predictable sets (see for example [2], or [28]) and pseudo-stopping times ([59]). **In this essay, we shall always note $L$ instead of $\rho$ for the end of an optional or predictable set[5] $\Gamma$**

$$L = \sup \left\{ t : \ (t, \omega) \in \Gamma \right\}.$$

Indeed, these random times, as will be clear in the sequel, have many interesting properties on their own, and hence, noting them differently will avoid confusion. Pseudo-stopping times have been discovered only recently: they have been introduced in [59], following Williams [74]:

**Definition 8.2.** We say that $\rho$ is a $(\mathcal{F}_t)$ pseudo-stopping time if for every $(\mathcal{F}_t)$-martingale $(M_t)$ in $\mathcal{H}^1$, we have

$$\mathbb{E} M_\rho = \mathbb{E} M_0. \tag{8.1}$$

**Remark 8.3.** It is equivalent to assume that (8.1) holds for bounded martingales, since these are dense in $\mathcal{H}^1$. It can also be proved that then (8.1) also holds for all uniformly integrable martingales (see [59]).

We shall in the sequel give a characterization of pseudo-stopping times but for now we indicate immediately that a class of pseudo-stopping times with respect to a filtration $(\mathcal{F}_t)$ which are not in general $(\mathcal{F}_t)$ stopping times may be obtained by considering stopping times with respect to a larger filtration $(\mathcal{G}_t)$ such that $(\mathcal{F}_t)$ is immersed in $(\mathcal{G}_t)$, i.e: every $(\mathcal{F}_t)$ martingale is a $(\mathcal{G}_t)$ martingale. This situation is described in ([22]) and refered to there as the $(H)$ hypothesis (this situation is discussed in more detail in subsection 9.3) . Here is a well known example: let $B_t = \left( B_t^1, \ldots, B_t^d \right)$ be a $d$-dimensional Brownian motion, and $R_t = |B_t|$, $t \geq 0$, its radial part; it is well known that

$$\left( \mathcal{R}_t \equiv \sigma \left\{ R_s, \ s \leq t \right\}, \ t \geq 0 \right),$$

the natural filtration of $R$, is immersed in $\left( \mathcal{B}_t \equiv \sigma \left\{ B_s, \ s \leq t \right\}, \ t \geq 0 \right)$, the natural filtration of $B$. Thus an example of $(\mathcal{R}_t)$ pseudo-stopping time is:

$$T_a^{(1)} = \inf \left\{ t : \ B_t^1 > a \right\}.$$

Recently, D. Williams [74] showed that with respect to the filtration $(\mathcal{F}_t)$ generated by a one dimensional Brownian motion $(B_t)_{t \geq 0}$, there exist pseudo-stopping times $\rho$ which are not $(\mathcal{F}_t)$ stopping times. D. Williams' example is the following: let

$$T_1 = \inf \left\{ t : \ B_t = 1 \right\}, \text{ and } \sigma = \sup \left\{ t < T_1 : \ B_t = 0 \right\};$$

then

$$\rho = \sup \left\{ s < \sigma : \ B_s = S_s \right\}, \text{ where } S_s = \sup_{u \leq s} B_u$$

is a pseudo-stopping time.

Now, we give the definitions of some sigma fields associated with arbitrary random times, following Chung and Doob ([23]):

---

[5]We can make the convention that $\sup \emptyset = -\infty$.

**Definition 8.4.** Three classical $\sigma$-fields associated with a filtration $(\mathcal{F}_t)$ and any random time $\rho$ are:

$$\begin{cases} \mathcal{F}_{\rho+} & = & \sigma\left\{z_\rho, \ (z_t) \text{ any } (\mathcal{F}_t) \text{ progressively measurable process}\right\}; \\ \mathcal{F}_{\rho} & = & \sigma\left\{z_\rho, \ (z_t) \text{ any } (\mathcal{F}_t) \text{ optional process}\right\}; \\ \mathcal{F}_{\rho-} & = & \sigma\left\{z_\rho, \ (z_t) \text{ any } (\mathcal{F}_t) \text{ predictable process}\right\}. \end{cases}$$

**Remark 8.5.** As usual, when dealing with predictable processes on $[0, \infty]$, we assume that there is a sigma field $\mathcal{F}_{0-} = \mathcal{F}_0$ in the filtration $(\mathcal{F}_t)$.

**Remark 8.6.** When $\rho$ is a stopping time, we have $\mathcal{F}_\rho = \mathcal{F}_{\rho+}$ and the definitions of $\mathcal{F}_\rho$ and $\mathcal{F}_{\rho-}$ coincide of course with the usual definitions of the sigma fields associated with a stopping time.

**Remark 8.7.** In general, $\mathcal{F}_{\rho+} \neq \mathcal{F}_\rho$; for example, take $\rho$ to be the last time before 1 a standard Brownian motion is equal to zero; then the sign of the excursion between $\rho$ and 1 is $\mathcal{F}_{\rho+}$ measurable and orthogonal to $\mathcal{F}_\rho$ (see [69], [80] or [28]). We shall have a nice discussion about the differences between these two sigma fields, related to Brownian filtrations, at the end of this section.

One must be very careful when comparing two such sigma fields. For example, $\rho \leq \rho'$ does not necessarily imply that $\mathcal{F}_\rho \subset \mathcal{F}_{\rho'}$. However, we have the following useful result:

**Theorem 8.8** ([28], p. 142)**.** *Let $\rho$ and $\rho'$ be two random times, such that $\rho \leq \rho'$. If $\rho$ is measurable with respect to $\mathcal{F}_{\rho'}$ (resp. $\mathcal{F}_{\rho'-}$), then*

$$\mathcal{F}_\rho \subset \mathcal{F}_{\rho'} \quad (\text{resp. } \mathcal{F}_{\rho-} \subset \mathcal{F}_{\rho'-}).$$

*The previous assumption is always satisfied if $\rho$ is the end of an optional (resp. predictable) set.*

When dealing with arbitrary random times, one often works under the following conditions:

- Assumption (**C**): all $(\mathcal{F}_t)$-martingales are <u>c</u>ontinuous (e.g: the Brownian filtration).
- Assumption (**A**): the random time $\rho$ <u>a</u>voids every $(\mathcal{F}_t)$-stopping time $T$, i.e. $\mathbb{P}\left[\rho = T\right] = 0$.

When we refer to assumptions (**CA**), this will mean that both the conditions (**C**) and (**A**) hold.

**Lemma 8.9.** *Under the condition* (**A**), *we have*

$$\mathcal{F}_\rho = \mathcal{F}_{\rho-}.$$

There is also another important family of random times, called honest times, and which in fact coincides with ends of optional sets.

**Definition 8.10.** Let $L$ be a random variable with values in $[0, \infty]$. $L$ is said to be honest if for every $t$, there exists an $(\mathcal{F}_t)$ measurable random variable $L_t$, such that $L = L_t$ on the set $\{L < t\}$.

Every stopping time is an honest time (take $L_t \equiv L \wedge t$). There are also examples of honest times which are not stopping times. For example, let $X$ be an adapted and continuous process and set: $\overline{X}_t = \sup_{s \leq t} X_s$, $\overline{X} = \sup_{s \geq 0} X_s$. Then the random variable

$$L = \inf \left\{ s : X_s = \overline{X} \right\}$$

is honest. Indeed, on the set $\{L < t\}$, we have $L = \inf \left\{ s : X_s = \overline{X}_t \right\}$, which is $(\mathcal{F}_t)$ measurable[6]. Now, we characterize honest times (see [41], [28], p. 137, or [68] p. 373):

**Theorem 8.11** ([25])**.** *Let $L$ be a random time. Then $L$ is an honest time if and only if there exists an optional set $H$ in $[0, \infty] \times \Omega$ such that $L$ is the end of $H$.*

**Remark 8.12.** Ends of optional sets $H \subset [0, \infty[ \times \Omega$ do not allow to construct all honest times (Jeulin [42], p.74).

We postpone examples to the next section where we are able to give more details.

### 8.2. Azéma's supermartingales and dual projections associated with random times

A few processes play a crucial role in the study of arbitrary random times:

- the $(\mathcal{F}_t)$ supermartingale

$$Z_t^\rho = \mathbb{P} \left[ \rho > t \mid \mathcal{F}_t \right] \tag{8.2}$$

chosen to be càdlàg, associated to $\rho$ by Azéma ([2]);
- the $(\mathcal{F}_t)$ dual optional and predictable projections of the process $1_{\{\rho \leq t\}}$, denoted respectively by $A_t^\rho$ and $a_t^\rho$;
- the càdlàg martingale

$$\mu_t^\rho = \mathbb{E} \left[ A_\infty^\rho \mid \mathcal{F}_t \right] = A_t^\rho + Z_t^\rho$$

which is in $\mathrm{BMO}(\mathcal{F}_t)$ (see [28] or [81]).

We also consider the Doob-Meyer decomposition of (8.2):

$$Z_t^\rho = m_t^\rho - a_t^\rho.$$

We can note that the supermartingale $(Z_t^\rho)$ is the optional projection of $\mathbf{1}_{[0,\rho[}$.

---

[6]We would like here to quote Paul André Meyer ([28], p.137):

Par exemple $X_t$ peut représenter le cours d'une certaine action à l'instant $t$, et $L$ est le moment idéal pour vendre son paquet d'actions. Tous les spéculateurs cherchent à connaître $L$ sans jamais y parvenir, d'où son nom de v.a. honnête.

**Lemma 8.13.** *If $\rho$ avoids any $(\mathcal{F}_t)$ stopping time (i.e. condition (**A**) is satisfied), then $A_t^\rho = a_t^\rho$ is continuous.*

Under condition (**C**), $A^\rho$ is predictable (recall that we proved that under condition (**C**) the predictable and optional sigma fields are equal) and consequently $A^\rho = a^\rho$.

Under conditions (**CA**), $Z^\rho$ is continuous.

First we give a result which is not so well known and which may turn out to be useful in modeling default times:

**Proposition 8.14** ([28], p.134)**.** *Let $\rho$ be an arbitrary random time. The sets $\{Z^\rho = 0\}$ and $\{Z_-^\rho = 0\}$ are both disjoint from the stochastic interval $[0, \rho[$, and have the same lower bound $T$, which is the smallest stopping time bigger than $\rho$.*

When studying random times which are not stopping times, one usually makes the assumption (**A**): consequently, **in the sequel, we make the assumption (A) unless stated otherwise.** The reader can refer to [28, 42] if he wants more general results.

### 8.2.1. The case of honest times

We now concentrate on honest times which are the best known random times after stopping times. We state a very important result of Azéma which has many applications and which is very useful in the theory of progressive enlargements of filtrations:

**Theorem 8.15** (Azéma [2])**.** *Let $L$ be an honest time that avoids $(\mathcal{F}_t)$ stopping times, and let*

$$Z_t^L = \mathbb{P}\left[L > t \mid \mathcal{F}_t\right].$$

*Let*

$$Z_t^L = \mu_t^L - A_t^L,$$

*denote its Doob-Meyer decomposition. Then $A_\infty$ follows the exponential law with parameter 1 and the measure $dA_t$ is carried by the set $\{t : Z_t = 1\}$. Moreover, $A$ does not increase after $L$, i.e. $A_L = A_\infty$.*

*Finally we have:*

$$L = \sup\left\{t : 1 - Z_t = 0\right\}.$$

Let us now give an example. Consider again a Bessel process of dimension $2(1 - \mu)$, starting from 0. Set

$$g_\mu \equiv \sup\left\{t \leq 1 : R_t = 0\right\},$$

and more generally, for $T > 0$, a fixed time,

$$g_\mu(T) \equiv \sup\left\{t \leq T : R_t = 0\right\}.$$

**Proposition 8.16.** *Let $\mu \in (0,1)$, and let $(R_t)$ be a Bessel process of dimension $\delta = 2(1-\mu)$. Then, we have:*

$$Z_t^{g_\mu} \equiv \mathbb{P}\left[g_\mu > t \mid \mathcal{F}_t\right] = \frac{1}{2^{\mu-1}\Gamma(\mu)} \int_{\frac{R_t}{\sqrt{1-t}}}^{\infty} dy\, y^{2\mu-1} \exp\left(-\frac{y^2}{2}\right).$$

*Proof.* We have:

$$Z_t^{g_\mu} = 1 - \mathbb{P}\left[g_\mu \le t \mid \mathcal{F}_t\right] = 1 - \mathbb{P}\left[d_t > 1 \mid \mathcal{F}_t\right],$$

where (using the Markov property),

$$d_t = \inf\{u \ge t;\ R_u = 0\} = t + \inf\{u \ge 0;\ R_{t+u} = 0\} = t + \widehat{H}_0,$$

with

$$\widehat{H}_0 \equiv \inf\left\{u \ge 0;\ \widehat{R}_u = 0;\ \widehat{R}_0 = R_t\right\},$$

i.e. $\widehat{H}_0$ is the first time when a Bessel process of dimension $\delta$, starting from $R_t$ (we call its law $\widehat{\mathbb{P}}$), hits 0. Thus, we have proved so far that:

$$Z_t^{g_\mu} = 1 - \widehat{\mathbb{P}}\left[\widehat{H}_0 > 1 - t\right]. \tag{8.3}$$

Now, following Borodin and Salminen ([21], p. 70-71), if for $-\nu > 0$, $\mathbb{P}_0^{(-\nu)}$ denotes the law of a Bessel process of parameter $-\nu$, starting from 0, then the law of $L_y \equiv \sup\{t:\ R_t = y\}$, is given by:

$$\mathbb{P}_0^{(-\nu)}\left(L_y \in dt\right) = \frac{y^{-2\nu}}{2^{-\nu}\Gamma(-\nu)\,t^{-\nu+1}} \exp\left(-\frac{y^2}{2t}\right) dt.$$

Now, from the time reversal property for Bessel processes ([21] p.70, or [69]), we have:

$$\widehat{\mathbb{P}}\left[\widehat{H}_0 \in dt\right] = \mathbb{P}_0^{(-\nu)}\left(L_{R_t} \in dt\right);$$

consequently, from (8.3), we have (recall $\mu = -\nu$):

$$Z_t^{g_\mu} = 1 - \frac{R_t^{2\mu}}{2^\mu \Gamma(\mu)} \int_{1-t}^{\infty} du\, \frac{\exp\left(-\frac{R_t^2}{2u}\right)}{u^{1+\mu}},$$

and the desired result is obtained by straightforward change of variables in the above integral. $\qquad\square$

**Remark 8.17.** The previous proof can be applied mutatis mutandis to obtain:

$$\mathbb{P}\left[g_\mu(T) > t \mid \mathcal{F}_t\right] = \frac{1}{2^{\mu-1}\Gamma(\mu)} \int_{\frac{R_t}{\sqrt{T-t}}}^{\infty} dy\, y^{2\mu-1} \exp\left(-\frac{y^2}{2}\right);$$

and

$$A_t^{g_\mu(T)} = \frac{1}{2^\mu \Gamma(1+\mu)} \int_0^{t \wedge T} \frac{dL_u}{(T-u)^\mu}.$$

**Remark 8.18.** It can be easily deduced from Proposition 8.16 that the dual predictable projection $A_t^{g_\mu}$ of $\mathbf{1}_{(g_\mu \leq t)}$ is:

$$A_t^{g_\mu} = \frac{1}{2^\mu \Gamma (1 + \mu)} \int_0^{t \wedge 1} \frac{dL_u}{(1 - u)^\mu}.$$

Indeed, it is a consequence of Itô's formula applied to $Z_t^{g_\mu}$ and the fact that $N_t \equiv R_t^{2\mu} - L_t$ is a martingale and $(dL_t)$ is carried by $\{t : R_t = 0\}$.

When $\mu = \frac{1}{2}$, $R_t$ can be viewed as $|B_t|$, the absolute value of a standard Brownian Motion. Thus, we recover as a particular case of our framework the celebrated example of the last zero before 1 of a standard Brownian Motion (see [42] p.124, or [81] for more references).

**Corollary 8.19.** *Let $(B_t)$ denote a standard Brownian Motion and let*

$$g \equiv \sup \{t \leq 1 : B_t = 0\}.$$

*Then:*

$$\mathbb{P} [g > t \mid \mathcal{F}_t] = \sqrt{\frac{2}{\pi}} \int_{\frac{|B_t|}{\sqrt{1-t}}}^{\infty} dy \exp \left( -\frac{y^2}{2} \right),$$

*and*

$$A_t^g = \sqrt{\frac{2}{\pi}} \int_0^{t \wedge 1} \frac{dL_u}{\sqrt{1 - u}}.$$

*Proof.* It suffices to take $\mu \equiv \frac{1}{2}$ in Proposition 8.16. $\qquad\square$

**Corollary 8.20.** *The variable*

$$\frac{1}{2^\mu \Gamma (1 + \mu)} \int_0^1 \frac{dL_u}{(1 - u)^\mu}$$

*is exponentially distributed with expectation 1; consequently, its law is independent of $\mu$.*

*Proof.* The random time $g_\mu$ is honest by definition (it is the end of a predictable set). It also avoids stopping times since $A_t^{g_\mu}$ is continuous (this can also be seen as a consequence of the strong Markov property for $R$ and the fact that 0 is instantaneously reflecting). Thus the result of the corollary is a consequence of Remark 8.18 following Proposition 8.16 and Lemma 8.15. $\qquad\square$

Given an honest time, it is not in general easy to compute its associated supermartingale $Z^L$. Hence it is important (in view of the theory of progressive enlargements of filtrations) to dispose some characterizations of Azéma's supermartingales which also provide a method way to compute them explicitly. We will give two results in this direction, borrowed from [63] and [60].

Let $(N_t)_{t \geq 0}$ be a continuous local martingale such that $N_0 = 1$, and $\lim_{t \to \infty} N_t = 0$. Let $S_t = \sup_{s \leq t} N_s$. We consider:

$$
\begin{aligned}
g &= \sup \{ t \geq 0 : \quad N_t = S_\infty \} \\
&= \sup \{ t \geq 0 : \quad S_t - N_t = 0 \} .
\end{aligned}
\tag{8.4}
$$

**Proposition 8.21** ([63])**.** *Consider the supermartingale*

$$
Z_t \equiv \mathbf{P} \left( g > t \mid \mathcal{F}_t \right) .
$$

*1. In our setting, the formula:*

$$
Z_t = \frac{N_t}{S_t}, \ t \geq 0
$$

*holds.*

*2. The Doob-Meyer additive decomposition of $(Z_t)$ is:*

$$
Z_t = \mathbf{E} \left[ \log S_\infty \mid \mathcal{F}_t \right] - \log \left( S_t \right) .
\tag{8.5}
$$

The above proposition gives a large family of examples. In fact, quite remarkably , every supermartingale associated with an honest time is of this form. More precisely:

**Theorem 8.22** ([63])**.** *Let $L$ be an honest time. Then, under the conditions* **(CA)**, *there exists a continuous and nonnegative local martingale $(N_t)_{t \geq 0}$, with $N_0 = 1$ and $\lim_{t \to \infty} N_t = 0$, such that:*

$$
Z_t = \mathbf{P} \left( L > t \mid \mathcal{F}_t \right) = \frac{N_t}{S_t}.
$$

We shall now outline a nontrivial consequence of Theorem 8.22 here. In [7], the authors are interested in giving explicit examples of dual predictable projections of processes of the form $\mathbf{1}_{L \leq t}$, where $L$ is an honest time. Indeed, these dual projections are natural examples of increasing injective processes (see [7] for more details and references). With Theorem 8.22, we have a complete characterization of such projections:

**Corollary 8.23.** *Assume the assumption* **(C)** *holds, and let $(C_t)$ be an increasing process. Then $C$ is the dual predictable projection of $\mathbf{1}_{g \leq t}$, for some honest time $g$ that avoids stopping times, if and only if there exists a continuous local martingale $N_t$ in the class $\mathcal{C}_0$ such that*

$$
C_t = \log S_t.
$$

Now let us give some examples.

**Example 8.24.** Let

$$
N_t \equiv B_t,
$$

where $(B_t)_{t \geq 0}$ is a Brownian Motion starting at 1, and stopped at $T_0 = \inf \{t : B_t = 0\}$. Let

$$S_t \equiv \sup_{s \leq t} B_s.$$

Let

$$g = \sup \{t : B_t = S_t\}.$$

Then

$$\mathbf{P}(g > t \mid \mathcal{F}_t) = \frac{B_t}{S_t}.$$

**Example 8.25.** Let

$$N_t \equiv \exp \left(2\nu B_t - 2\nu^2 t\right),$$

where $(B_t)$ is a standard Brownian Motion, and $\nu > 0$. We have:

$$S_t = \exp \left( \sup_{s \leq t} 2\nu \left(B_s - \nu s\right) \right),$$

and

$$g = \sup \left\{ t : (B_t - \nu t) = \sup_{s \geq 0} (B_s - \nu s) \right\}.$$

Consequently,

$$\mathbf{P}(g > t \mid \mathcal{F}_t) = \exp \left( 2\nu \left( (B_s - \nu s) - \sup_{s \leq t} (B_s - \nu s) \right) \right).$$

**Example 8.26.** Now, we consider $(R_t)$, a transient diffusion with values in $[0, \infty)$, which has $\{0\}$ as entrance boundary. Let $s$ be a scale function for $R$, which we can choose such that:

$$s(0) = -\infty, \text{ and } s(\infty) = 0.$$

Then, under the law $\mathbf{P}_x$, $x > 0$, the local martingale $\left( N_t = \dfrac{s(R_t)}{s(x)}, \ t \geq 0 \right)$ satisfies the required conditions of Proposition 8.21, and we have:

$$\mathbf{P}_x (g > t | \mathcal{F}_t) = \frac{s(R_t)}{s(I_t)}$$

where

$$g = \sup \{t : R_t = I_t\},$$

and

$$I_t = \inf_{s \leq t} R_s.$$

Theorem 8.22 is a multiplicative characterization; now we shall give an additive one.

**Theorem 8.27** ([60])**.** *Again, we assume that the conditions* (**CA**) *hold. Let* $(X_t)$ *be a submartingale of the class* $(\Sigma_c D)$ *satisfying:* $\lim_{t \to \infty} X_t = 1$. *Let*

$$L = \sup \{ t : X_t = 0 \}.$$

*Then* $(X_t)$ *is related to the Azéma's supermartingale associated with* $L$ *in the following way:*

$$X_t = 1 - Z_t^L = \mathbb{P}(L \leq t | \mathcal{F}_t).$$

*Consequently, if* $(Z_t)$ *is a nonnegative supermartingale, with* $Z_0 = 1$, *then,* $Z$ *may be represented as* $\mathbf{P}(L > t | \mathcal{F}_t)$, *for some honest time* $L$ *which avoids stopping times, if and only if* $(X_t \equiv 1 - Z_t)$ *is a submartingale of the class* $(\Sigma)$, *with the limit condition:*

$$\lim_{t \to \infty} X_t = 1.$$

Now, we give some fundamental examples:

**Example 8.28.** First, consider $(B_t)$, the standard Brownian Motion, and let $T_1 = \inf \{ t \geq 0 : B_t = 1 \}$. Let $\sigma = \sup \{ t < T_1 : B_t = 0 \}$. Then $B_{t \wedge T_1}^+$ satisfies the conditions of Theorem 8.27, and hence:

$$\mathbb{P}(\sigma \leq t | \mathcal{F}_t) = B_{t \wedge T_1}^+ = \int_0^{t \wedge T_1} \mathbf{1}_{B_u > 0} \, dB_u + \frac{1}{2} \ell_{t \wedge T_1},$$

where $(\ell_t)$ is the local time of $B$ at 0. This example plays an important role in the celebrated Williams' path decomposition for the standard Brownian Motion on $[0, T_1]$.

One can also consider $T_{\pm 1} = \inf \{ t \geq 0 : |B_t| = 1 \}$ and $\tau = \sup \{ t < T_{\pm 1} : |B_t| = 0 \}$. $|B_{t \wedge T_{\pm 1}}|$ satisfies the conditions of Theorem 8.27, and hence:

$$\mathbb{P}(\tau \leq t | \mathcal{F}_t) = |B_{t \wedge T_{\pm 1}}| = \int_0^{t \wedge T_{\pm 1}} sgn(B_u) \, dB_u + \ell_{t \wedge T_{\pm 1}}.$$

**Example 8.29.** Let $(Y_t)$ be a real continuous recurrent diffusion process, with $Y_0 = 0$. Then from the general theory of diffusion processes, there exists a unique continuous and strictly increasing function $s$, with $s(0) = 0$, $\lim_{x \to +\infty} s(x) = +\infty$, $\lim_{x \to -\infty} s(x) = -\infty$, such that $s(Y_t)$ is a continuous local martingale. Let

$$T_1 \equiv \inf \{ t \geq 0 : Y_t = 1 \}.$$

Now, if we define

$$X_t \equiv \frac{s(Y_{t \wedge T_1})^+}{s(1)},$$

we easily note that $X$ is a local submartingale of the class $(\Sigma_c)$ which satisfies the hypotheses of Theorem 8.27. Consequently, if we note

$$\sigma = \sup \{ t < T_1 : Y_t = 0 \},$$

we have:

$$\mathbb{P}\left(\sigma \le t | \mathcal{F}_t\right) = \frac{s\left(Y_{t \wedge T_1}\right)^+}{s\left(1\right)}.$$

**Example 8.30.** Now let $(M_t)$ be a positive local martingale, such that: $M_0 = x$, $x > 0$ and $\lim_{t \to \infty} M_t = 0$. Then, Tanaka's formula shows us that $\left(1 - \dfrac{M_t}{y} \wedge 1\right)$, for $0 \le y \le x$, is a local submartingale of the class $(\Sigma_c)$ satisfying the assumptions of Theorem 8.27, and hence with

$$g = \sup\left\{t : \ M_t = y\right\},$$

we have:

$$\mathbb{P}\left(g > t | \mathcal{F}_t\right) = \frac{M_t}{y} \wedge 1 = 1 + \frac{1}{y} \int_0^t \mathbf{1}_{(M_u < y)} \, dM_u - \frac{1}{2y} L_t^y,$$

where $(L_t^y)$ is the local time of $M$ at $y$.

**Example 8.31.** As an illustration of the previous example, consider $(R_t)$, a transient diffusion with values in $[0, \infty)$, which has $\{0\}$ as entrance boundary. Let $s$ be a scale function for $R$, which we can choose such that:

$$s\left(0\right) = -\infty, \ \text{and} \ s\left(\infty\right) = 0.$$

Then, under the law $\mathbb{P}_x$, for any $x > 0$, the local martingale $(M_t = -s\left(R_t\right))$ satisfies the conditions of the previous example and for $0 \le x \le y$, we have:

$$\mathbb{P}_x\left(g_y > t | \mathcal{F}_t\right) = \frac{s\left(R_t\right)}{s\left(y\right)} \wedge 1 = 1 + \frac{1}{s\left(y\right)} \int_0^t \mathbf{1}_{(R_u > y)} \, d\left(s\left(R_u\right)\right) + \frac{1}{2s\left(y\right)} L_t^{s(y)},$$

where $\left(L_t^{s(y)}\right)$ is the local time of $s\left(R\right)$ at $s\left(y\right)$, and where

$$g_y = \sup\left\{t : \ R_t = y\right\}.$$

This last formula was the key point for deriving the distribution of $g_y$ in [67], Theorem 6.1, p.326.

### 8.2.2. The case of pseudo-stopping times

In this paragraph, we give some characteristic properties and some examples of pseudo-stopping times. We do not assume here that condition (**A**) holds, but we assume that $\mathbb{P}\left[\rho = \infty\right] = 0$.

**Theorem 8.32** ([59]). *The following properties are equivalent:*

1. *$\rho$ is a $(\mathcal{F}_t)$ pseudo-stopping time, i.e (8.1) is satisfied;*
2. *$A_\infty^\rho \equiv 1$, a.s*

**Remark 8.33.** We shall give a more complete version of Theorem 8.32 in the section on progressive expansions of filtrations.

*Proof.* We have:

$$\mathbb{E}\left[M_\rho\right] = \mathbb{E}\left[\int_0^\infty M_s dA_s^\rho\right] = \mathbb{E}\left[M_\infty A_\infty^\rho\right].$$

Hence,

$$\mathbb{E}\left[M_\rho\right] = \mathbb{E}\left[M_\infty\right] \Leftrightarrow \mathbb{E}\left[M_\infty\left(A_\infty^\rho - 1\right)\right] = 0,$$

and the announced equivalence follows now easily. □

**Remark 8.34.** More generally, the approach adopted in the proof can be used to solve the equation

$$\mathbb{E}\left[M_\rho\right] = \mathbb{E}\left[M_\infty\right],$$

where the random time $\rho$ is fixed and where the unknown are martingales in $\mathcal{H}^1$. For more details and resolutions of such equations, see [64].

**Corollary 8.35.** *Under the assumptions of Theorem 8.32, $Z_t^\rho = 1 - A_t^\rho$ is a decreasing process. Furthermore, if $\rho$ avoids stopping times, then $(Z_t^\rho)$ is continuous.*

*Proof.* The follows from the fact that $\mu_t^\rho = \mathbb{E}\left[A_\infty^\rho | \mathcal{F}_t\right] = 1$. □

**Remark 8.36.** In fact, we shall see in next section, that under condition (**C**), $\rho$ is a pseudo-stopping time if and only if $(Z_t^\rho)$ is a predictable decreasing process.

For honest times, Azéma proved that $A_L$ follows the standard exponential law. For pseudo-stopping times, we have:

**Proposition 8.37** ([59])**.** *For simplicity, we shall write $(Z_u)$ instead of $(Z_u^\rho)$. Under condition (A), for all bounded $(\mathcal{F}_t)$ martingales $(M_t)$, and all bounded Borel measurable functions $f$, one has:*

$$\begin{aligned}
\mathbb{E}\left[M_\rho f\left(Z_\rho\right)\right] &=& \mathbb{E}\left[M_0\right] \int_0^1 f\left(x\right) dx \\
&=& \mathbb{E}\left[M_\rho\right] \int_0^1 f\left(x\right) dx.
\end{aligned}$$

*Consequently, $Z_\rho$ follows the uniform law on $(0,1)$.*

*Proof.* Under our assumptions, we have

$$
\begin{aligned}
\mathbb{E}\left[M_\rho f\left(Z_\rho\right)\right] &= \mathbb{E}\left[\int_0^\infty M_u f\left(Z_u\right) dA_u^\rho\right] \\
&= \mathbb{E}\left[\int_0^\infty M_u f\left(1 - A_u^\rho\right) dA_u^\rho\right] \\
&= \mathbb{E}\left[M_\infty \int_0^\infty f\left(1 - A_u^\rho\right) dA_u^\rho\right] \\
&= \mathbb{E}\left[M_\infty \int_0^1 f\left(1 - x\right) dx\right] \\
&= \mathbb{E}\left[M_\infty \int_0^1 f\left(x\right) dx\right].
\end{aligned}
$$

$\square$

Now, we give a systematic construction for pseudo-stopping times, generalizing D. Williams's example. We assume we are given an honest time $L$ and that conditions (**CA**) hold (the condition (**A**) holds with respect to $L$). Then the following holds:

**Proposition 8.38** ([59]). *(i) $I_L = \inf_{u \leq L} Z_u^L$ is uniformly distributed on $[0, 1]$;*

*(ii) The supermartingale $Z_t^\rho = P\left[\rho > t \mid \mathcal{F}_t\right]$ associated with $\rho$ is given by*

$$
Z_t^\rho = \inf_{u \leq t} Z_u^L.
$$

*As a consequence, $\rho$ is a $(\mathcal{F}_t)$ pseudo-stopping time.*

*Proof.* For simplicity, we write $Z_t$ for $Z_t^L$. (i) Let

$$
T_b = \inf\left\{t, \quad Z_t \leq b\right\}, \qquad 0 < b < 1,
$$

then

$$
\mathbb{P}\left[I_L \leq b\right] = \mathbb{P}\left[T_b < L\right] = \mathbb{E}\left[Z_{T_b}\right] = b.
$$

(*ii*) Note that for every $(\mathcal{F}_t)$ stopping time $T$, we have

$$
\{T < \rho\} = \left\{T' < L\right\}
$$

where

$$
T' = \inf\left\{t > T, \quad Z_t \leq \inf_{s \leq T} Z_s\right\}.
$$

Consequently, we have

$$
\mathbb{E}\left[Z_T^\rho\right] = \mathbb{P}\left[T < \rho\right] = \mathbb{P}\left[T' < L\right] = \mathbb{E}\left[Z_{T'}\right] = \mathbb{E}\left[\inf_{u \leq T} Z_u\right],
$$

which yields:

$$\mathbb{E}\left[Z_T^\rho \mathbf{1}_{\{T<\infty\}}\right] = \mathbb{E}\left[\inf_{u\leq T} Z_u \mathbf{1}_{\{T<\infty\}}\right],$$

since $(Z_u^\rho)$ and $(Z_u)$ converge to 0 as $u \to \infty$. We now deduce the desired result from the optional section theorem. $\qquad\square$

Now, let us give some examples of pseudo stopping times. We shall use the supermartingales associated with honest times we have computed in the previous paragraph, and for simplicity, we write $Z_t$ for $Z_t^L$.

1. First let us check that we recover the example of D. Williams from the proposition. With the notations of D. Williams's example $(L = \sigma)$, it is not hard to see that (see [70])

$$Z_t = 1 - B_{t \wedge T_1}^+.$$

   Hence

$$\rho = \sup\left\{s < \sigma : \ B_s = S_s\right\}.$$

2. Consider $(R_t)_{t\geq 0}$ a three dimensional Bessel process, starting from zero, its filtration $(\mathcal{F}_t)$, and

$$L = L_1 = \sup\left\{t : \ R_t = 1\right\}.$$

   Then

$$\rho = \sup\left\{t < L : \ R_t = \sup_{u\leq L} R_u\right\}, \qquad (8.6)$$

   is a $(\mathcal{F}_t)$ pseudo-stopping time. This follows from the fact that

$$Z_t^L = 1 \wedge \frac{1}{R_t},$$

   hence (8.6) is equivalent to:

$$\rho = \sup\left\{t < L : \ Z_t^L = \inf_{u\leq L} Z_u^L\right\},$$

   and from the above proposition:

$$Z_t^\rho = 1 \wedge \left(\frac{1}{\sup\limits_{u\leq t} R_u}\right).$$

3. Similarly, with our previous notations on Bessel processes of dimension $2(1-\mu)$, $\mu \in (0,1)$, define:

$$\rho \equiv \sup\left\{t < g_\mu : \ \frac{R_t}{\sqrt{1-t}} = \sup_{u<g_\mu} \frac{R_u}{\sqrt{1-u}}\right\}.$$

   Then, $\rho$ is a pseudo-stopping time.

### 8.3. Honest times and Strong Brownian Filtrations

In this paragraph, we shall describe a very nice and very difficult recent result on Strong Brownian Filtrations. We will not go into details; the aim here is just to show a powerful application of non stopping times. We give the references for further details.

**Definition 8.39** (Strong Brownian Filtration). A Strong Brownian Filtration is a filtration which is generated by a standard Brownian Motion.

**Definition 8.40** (Weak Brownian Filtration). A filtration $(\mathcal{F}_t)$ is a Weak Brownian Filtration if there exists an $(\mathcal{F}_t)$ Brownian Motion $(\beta_t)$ such that every $(\mathcal{F}_t)$ local martingale $(M_t)$ may be represented as

$$M_t = c + \int_0^t m_s d\beta_s, \quad t \geq 0$$

for some $(\mathcal{F}_t)$ predictable process $(m_t)$.

Of course, any Strong Brownian Filtration is a Weak Brownian Filtration. The converse is not true. An example of a Weak Brownian Filtration which is not a Strong Brownian Filtration may be obtained considering Walsh's Brownian Motion. We introduce informally Walsh's Brownian Motion. Walsh's Brownian Motion $(Z_t)$ is a Feller process taking values in $N$ (half-lines) rays $(I_i;\ i = 1, \ldots, N)$ of the plane, all meeting at 0. Informally, $(Z_t)$ behaves like a Brownian Motion when it is away from 0, and when it reaches 0, its chooses its ray with equal probability $1/N$ (more generally it chooses its $i^{\text{th}}$ ray $I_i$ with with probability $p_i > 0$, and $\sum_{i=1}^n p_i = 1$). This description is not rigorous, since 0 is regular for itself (with respect to the Markov process $(Z_t)$), but it may be made rigorous using excursion theory (see [15]).

Moreover, it is shown in [15] that all martingales with respect to the natural filtration $(\mathcal{F}_t^Z)$ of $Z$ may be represented as stochastic integrals with respect to the Brownian Motion:

$$W_t = |Z_t| - \frac{1}{2}L_t,$$

where $(L_t)$ is the local time at 0 of the reflecting Brownian Motion $|Z|$.

It had been an open question for a long time to know whether or not $(\mathcal{F}_t^Z)$ is a Strong Brownian Filtration. The answer was given by Tsirelson in 1997 ([72]):

**Theorem 8.41.** *For $N \geq 3$, $(\mathcal{F}_t^Z)$ is not a Strong Brownian Filtration.*

Another proof of this theorem was given in 1998 by Barlow, Emery, Knight, Song and Yor ([18]). We give the logic of their proof. First, it was shown by Barlow, Pitman and Yor ([15]) that with

$$\gamma \equiv \sup\{t < 1:\ Z_t = 0\},$$

we have

$$\mathcal{F}_{\gamma+} = \mathcal{F}_\gamma \bigvee \sigma\{(Z_i)\,;\ i = 1, \ldots N\}.$$

To obtain Theorem 8.41, Barlow et alii. proved in [18] a result conjectured earlier by Barlow:

**Theorem 8.42** ([18]). *If $(\mathcal{F}_t)$ is a Strong Brownian Filtration and $L$ the end of a predictable set, then:*

$$\mathcal{F}_{L+} = \mathcal{F}_L \bigvee \sigma(A)$$

*with at most one non trivial set $A \in \mathcal{F}_{L+}$.*

## 9. The enlargements of filtrations

The aim of this section is to present the theory of enlargements of filtrations, and to give some applications. The main question is the following: how are semimartingales modified when considered as stochastic processes in a larger filtration $(\mathcal{G}_t)$ than the initial one $(\mathcal{F}_t)$ (i.e. for all $t \geq 0$, $\mathcal{F}_t \subset \mathcal{G}_t$)? A first result in this direction (in fact in the reverse direction) is a theorem of Stricker which we shall sometimes use in the sequel:

**Theorem 9.1** (Stricker [71]). *Let $(\mathcal{F}_t)$ and $(\mathcal{G}_t)$ be two filtrations such that for all $t \geq 0$, $\mathcal{F}_t \subset \mathcal{G}_t$. If $(X_t)$ is a $(\mathcal{G}_t)$ semimartingale which is $(\mathcal{F}_t)$ adapted, then it is also an $(\mathcal{F}_t)$ semimartingale.*

Given a filtered probability space $(\Omega, \mathcal{F}, (\mathcal{F}_t), \mathbb{P})$, there are essentially two ways of enlarging filtrations:

- *initial enlargements*, for which $\mathcal{G}_t = \mathcal{F}_t \bigvee \mathcal{H}$, i.e. the new information $\mathcal{H}$ is brought in at the origin of time; and
- *progressive enlargements*, for which $\mathcal{G}_t = \mathcal{F}_t \bigvee \mathcal{H}_t$, i.e. the new information is brought in progressively as the time $t$ increases.

We shall try to characterize situations when every $(\mathcal{F}_t)$ semimartingale $X$ remains a $(\mathcal{G}_t)$ semimartingale and then find the decomposition of $X$ as a $(\mathcal{G}_t)$ semimartingale. This situation is described as the $(H')$ hypothesis:

**Definition 9.2.** We shall say that the pair of filtrations $(\mathcal{F}_t, \mathcal{G}_t)$ satisfies the $(H')$ hypothesis if every $(\mathcal{F}_t)$ (semi)martingale is a $(\mathcal{G}_t)$ semimartingale.

**Remark 9.3.** In fact it suffices to check that every $(\mathcal{F}_t)$ martingale is a $(\mathcal{G}_t)$ semimartingale.

When the $(H')$ hypothesis is not satisfied, we shall try to find some conditions under which an $(\mathcal{F}_t)$ martingale is a $(\mathcal{G}_t)$ semimartingale.

Of course, the problem does not have a solution in the generality presented above. For the initial enlargement case, we shall deal with the case when $\mathcal{H}$ is the sigma field generated by a random variable $Z$ and for the progressive enlargement case, we shall take $\mathcal{H}_t = \sigma\{\rho \wedge t\}$, where $\rho$ is a random time, so that $(\mathcal{G}_t)$ is the smallest filtration which contains $(\mathcal{F}_t)$ and which makes $\rho$ a stopping time. All the results in the sequel originate from the works of Barlow, Jeulin, Jacod and Yor (see [42] or [45] for a complete account and more references; see also [81, 68]).

### *9.1. Initial enlargements of filtrations*

The theory of initial enlargements of filtrations is better known than the progressive enlargements of filtrations. The main results can be found in [42], [45], [81] or [68]. We first give a theoretical result of Jacod with some applications, and then we give a general method for obtaining the decomposition of a local martingale in the enlarged filtrations in some situations where Jacod's results do not apply.

Let $(\Omega, \mathcal{F}, (\mathcal{F}_t), \mathbb{P})$ be a filtered probability space satisfying the usual assumptions. Let $Z$ be an $\mathcal{F}$ measurable random variable. Define

$$\mathcal{G}_t = \cap_{\varepsilon > 0} \left( \mathcal{F}_{t+\varepsilon} \bigvee \sigma \{Z\} \right).$$

The conditional laws of $Z$ given $\mathcal{F}_t$, for $t \geq 0$ play a crucial role in initial enlargements.

**Theorem 9.4** (Jacod's criterion)**.** *Let $Z$ be an $\mathcal{F}$ measurable random variable and let $Q_t (\omega, dx)$ denote the regular conditional distribution of $Z$ given $\mathcal{F}_t$, $t \geq 0$. Suppose that for each $t \geq 0$, there exists a positive $\sigma$-finite measure $\eta_t (dx)$ (on $(\mathbb{R}, \mathcal{B}(\mathbb{R}))$) such that*

$$Q_t (\omega, dx) \ll \eta_t (dx) \text{ a.s.}$$

*Then every $(\mathcal{F}_t)$ semimartingale is a $(\mathcal{G}_t)$ semimartingale.*

*Proof.* See [38] (or [68] for an English reference). $\qquad\square$

**Remark 9.5.** In fact this theorem still holds for random variables with values in a standard Borel space. Moreover, the existence of the $\sigma$-finite measure $\eta_t (dx)$ is equivalent to the existence of one positive $\sigma$-finite measure $\eta (dx)$ such that $Q_t (\omega, dx) \ll \eta (dx)$ and in this case $\eta$ can be taken to be the distribution of $Z$.

Now we give classical corollaries of Jacod's theorem.

**Corollary 9.6.** *Let $Z$ be independent of $\mathcal{F}_\infty$. Then every $(\mathcal{F}_t)$ semimartingale is a $(\mathcal{G}_t)$ semimartingale.*

*Proof.* It suffices (with the notations of Theorem 9.4) to note that $Q_t (\omega, dx) = \eta (dx)$, where $\eta (dx)$ is the law of $Z$. $\qquad\square$

**Corollary 9.7.** *Let $Z$ be a random variable taking on only a countable number of values. Then every $(\mathcal{F}_t)$ semimartingale is a $(\mathcal{G}_t)$ semimartingale.*

*Proof.* If we note

$$\eta (dx) = \sum_{k=1}^{\infty} \mathbb{P} (Z = x_k) \, \delta_{x_k} (dx),$$

where $\delta_{x_k} (dx)$ is the Dirac measure at $x_k$, the law of $Z$, then $Q_t (\omega, dx)$ is absolutely continuous with respect to $\eta$ with Radon-Nikodym density:

$$\sum_{k=1}^{\infty} \frac{\mathbb{P} (Z = x_k | \mathcal{F}_t)}{\mathbb{P} (Z = x_k)} \mathbf{1}_{x = x_k}.$$

Now the result follows from Jacod's theorem. □

Now we give a theorem dealing with initial enlargement with $A_\infty^L$, when $L$ is an honest time. This theorem is important for the following reasons:

- The results of Jacod do not apply;
- In fact, the method we shall use applies to many other situations, where the theorems of Jacod do not apply.

Let $(\Omega, \mathcal{F}, (\mathcal{F}_t), \mathbb{P})$ be a filtered probability space satisfying the usual assumptions and let $L$ be an honest time. We assume for simplicity that the conditions (**CA**) hold. In the sequel, we shall note $A_t, Z_t, \mu_t$ for $A_t^L, Z_t^L, \mu_t^L$. Let us define the new filtration

$$\mathcal{F}_t^{\sigma(A_\infty)} \equiv \bigcap_{\varepsilon > 0} \left( \mathcal{F}_{t+\varepsilon} \vee \sigma\left( A_\infty \right) \right),$$

which satisfies the usual assumptions.

We first need the conditional laws of $A_\infty$ which were obtained under conditions (**CA**) in [6] and in a more general setting and by different methods in [61].

**Proposition 9.8** ([61],[6]). *Let $G$ be a Borel bounded function. Define:*

$$M_t^G \equiv \mathbf{E}\left( G\left( A_\infty \right) | \mathcal{F}_t \right).$$

*Then,*

$$M_t^G = F\left( A_t \right) - \left( F\left( A_t \right) - G\left( A_t \right) \right)\left( 1 - Z_t \right),$$

*where*

$$F\left( x \right) = \exp\left( x \right) \int_x^\infty dy \exp\left( -y \right) G\left( y \right).$$

*Moreover, $\left( M_t^G \right)$ has the following stochastic integral representation:*

$$M_t^G = \mathbf{E}\left[ G\left( A_\infty \right) \right] + \int_0^t \left( F - G \right)\left( A_u \right) d\mu_u.$$

Now, define, for $G$ any Borel bounded function,

$$\lambda_t\left( G \right) \equiv M_t^G = F\left( A_t \right) - \left( F\left( A_t \right) - G\left( A_t \right) \right)\left( 1 - Z_t^g \right).$$

From Proposition 9.8, we also have:

$$
\begin{aligned}
\lambda_t\left( G \right) &= \mathbf{E}\left[ G\left( A_\infty \right) \right] + \int_0^t \left( F - G \right)\left( A_s \right) d\mu_s \\
&\equiv \mathbf{E}\left[ G\left( A_\infty \right) \right] + \int_0^t \dot{\lambda}_s\left( G \right) d\mu_s.
\end{aligned}
$$

Hence we have:

$$\lambda_t\left( G \right) = \int \lambda_t\left( dx \right) G\left( x \right),$$

with

$$\lambda_t (dx) = (1 - Z_t) \, \delta_{A_t} (dx) + Z_t \exp (A_t) \, \mathbf{1}_{(A_t, \infty)} (x) \exp (-x) \, dx,$$

where $\delta_{A_t}$ denotes the Dirac mass at $A_t$. Similarly, we have:

$$\dot{\lambda}_t (G) = \int \dot{\lambda}_t (dx) \, G (x),$$

with:

$$\dot{\lambda}_t (dx) = -\delta_{A_t} (dx) + \exp (A_t) \, \mathbf{1}_{(A_t, \infty)} (x) \exp (-x) \, dx.$$

It then follows that:

$$\dot{\lambda}_t (dx) = \lambda_t (dx) \, \rho (x, t), \tag{9.1}$$

with

$$\rho (x, t) = \frac{1}{Z_t} \mathbf{1}_{\{x > A_t\}} - \frac{1}{1 - Z_t} \mathbf{1}_{\{x = A_t\}}. \tag{9.2}$$

Now we can state our result about initial expansion with $A_\infty$, which was first obtained by Jeulin ([42]), but the proof we shall present is borrowed from [60].

**Theorem 9.9.** *Let $L$ be an honest time. We assume, as usual, that the conditions **(CA)** hold. Then, every $(\mathcal{F}_t)$ local martingale $M$ is an $\left( \mathcal{F}_t^{\sigma(A_\infty)} \right)$ semimartingale and decomposes as:*

$$M_t = \widetilde{M}_t + \int_0^t \mathbf{1}_{\{L > s\}} \frac{d\langle M, \mu \rangle_s}{Z_s} - \int_0^t \mathbf{1}_{\{L \le s\}} \frac{d\langle M, \mu \rangle_s}{1 - Z_s}, \tag{9.3}$$

*where $\left( \widetilde{M}_t \right)_{t \ge 0}$ denotes an $\left( \mathcal{F}_t^{\sigma(A_\infty)} \right)$ local martingale.*

*Proof.* We can first assume that $M$ is an $L^2$ martingale; the general case follows by localization. Let $\Lambda_s$ be an $\mathcal{F}_s$ measurable set, and take $t > s$. Then, for any bounded test function $G$, we have:

$$\mathbf{E} \left( \mathbf{1}_{\Lambda_s} G \left( A_\infty \right) \left( M_t - M_s \right) \right) = \mathbf{E} \left( \mathbf{1}_{\Lambda_s} \left( \lambda_t (G) \, M_t - \lambda_s (G) \, M_s \right) \right) \tag{9.4}$$

$$= \mathbf{E} \left( \mathbf{1}_{\Lambda_s} \left( \langle \lambda (G), M \rangle_t - \langle \lambda (G), M \rangle_s \right) \right) \tag{9.5}$$

$$= \mathbf{E} \left( \mathbf{1}_{\Lambda_s} \left( \int_s^t \dot{\lambda}_u (G) \, d\langle M, \mu \rangle_u \right) \right) \tag{9.6}$$

$$= \mathbf{E} \left( \mathbf{1}_{\Lambda_s} \left( \int_s^t \int \lambda_u (dx) \, \rho (x, u) \, G (x) \, d\langle M, \mu \rangle_u \right) \right)$$

$$= \mathbf{E} \left( \mathbf{1}_{\Lambda_s} \left( \int_s^t d\langle M, \mu \rangle_u \rho \left( A_\infty, u \right) \right) \right). \tag{9.7}$$

But from (9.2), we have:

$$\rho \left( A_\infty, t \right) = \frac{1}{Z_t} \mathbf{1}_{\{A_\infty > A_t\}} - \frac{1}{1 - Z_t} \mathbf{1}_{\{A_\infty = A_t\}}.$$

It now suffices to notice that $(A_t)$ is constant after $L$ and $L$ is the first time when $A_\infty = A_t$, or in other words (for example, see [28] p. 134):

$$\mathbf{1}_{\{A_\infty > A_t\}} = \mathbf{1}_{\{L > t\}}, \text{ and } \mathbf{1}_{\{A_\infty = A_t\}} = \mathbf{1}_{\{L \leq t\}}.$$

$\square$

Let us emphasize again that the method we have used here applies to many other situations, where the theorems of Jacod do not apply. Each time the different relationships we have just mentioned between the quantities: $\lambda_t(G)$, $\dot{\lambda}_t(G)$, and $\lambda_t(dx)$, $\dot{\lambda}_t(dx)$, $\rho(x,t)$, hold, the above method and decomposition formula apply. Moreover, the condition **(C)** can be dropped and it is enough to have only a stochastic integral representation for $\lambda_t(G)$ (see [63] for a discussion). In the case of enlargement with $A_\infty$, everything is nice since every $(\mathcal{F}_t)$ local martingale $M$ is an $\left(\mathcal{F}_t^{\sigma(A_\infty)}\right)$ semimartingale. Sometimes, an integrability condition is needed as is shown by the following example.

**Example 9.10** ([81], p.34)**.** Let $Z = \int_0^\infty \varphi(s) \, dB_s$, for some $\varphi \in L^2(\mathbb{R}_+, ds)$. Recall that

$$\mathcal{G}_t = \cap_{\varepsilon > 0} \left( \mathcal{F}_{t+\varepsilon} \bigvee \sigma\{Z\} \right).$$

We wish to address the following question: is $(B_t)$ a $(\mathcal{G}_t)$ semimartingale?

The above method applies step by step: it is easy to compute $\lambda_t(dx)$, since conditionally on $\mathcal{F}_t$, $Z$ is gaussian, with mean $m_t = \int_0^t \varphi(s) \, dB_s$, and variance $\sigma_t^2 = \int_0^t \varphi^2(s) \, ds$. Consequently, the absolute continuity requirement (9.1) is satisfied with:

$$\rho(x,t) = \varphi(s) \frac{x - m_s}{\sigma_s^2}.$$

But here, the arguments in the proof of Theorem 9.9 (replace $M$ with $B$) do not always work since the quantities involved there (equations (9.4) to (9.7)) might be infinite; hence we have to impose an integrability condition. For example, if we assume that

$$\int_0^t \frac{|\varphi(s)|}{\sigma_s} ds < \infty,$$

then $(B_t)$, is a $(\mathcal{G}_t)$ semimartingale with canonical decomposition:

$$B_t = B_0 + \widetilde{B}_t + \int_0^t ds \frac{\varphi(s)}{\sigma_s^2} \left( \int_s^\infty \varphi(u) \, dB_u \right),$$

where $\left(\widetilde{B}_t\right)$ is a $(\mathcal{G}_t)$ Brownian Motion.

As a particular case, we may take: $Z = B_{t_0}$, for some fixed $t_0$. The above formula then becomes:

$$B_t = B_0 + \widetilde{B}_t + \int_0^{t \wedge t_0} ds \frac{B_{t_0} - B_s}{t_0 - s},$$

where $\left(\widetilde{B}_t\right)$ is a $(\mathcal{G}_t)$ Brownian Motion. In particular, $\left(\widetilde{B}_t\right)$ is independent of $\mathcal{G}_0 = \sigma\{B_{t_0}\}$, so that conditionally on $B_{t_0} = y$, or equivalently, when $(B_t,\ t \leq t_0)$ is considered under the bridge law $\mathbf{P}^{t_0}_{x,y}$, its canonical decomposition is:

$$B_t = x + \widetilde{B}_t + \int_0^t ds \frac{y - B_s}{t_0 - s},$$

where $\left(\widetilde{B}_t,\ t \leq t_0\right)$ is now a $\left(\mathbf{P}^{t_0}_{x,y}; (\mathcal{F}_t)\right)$ Brownian Motion.

**Example 9.11.** For more examples of initial enlargements using this method, see the forthcoming book [50].

### 9.2. Progressive enlargements of filtrations

The theory of progressive enlargements of filtrations was originally motivated by a paper of Millar [57] on random times and decomposition theorems. It was first independently developed by Barlow [13] and Yor [77], and further developed by Jeulin and Yor [43] and Jeulin [41, 42]. For further developments and details, the reader can also refer to [45] which is written in French or to [81, 50] or [68] chapter VI. for an English text.

Let $(\Omega, \mathcal{F}, (\mathcal{F}_t), \mathbb{P})$ be a filtered probability space satisfying the usual assumptions, and for simplicity (and because it is always the case with practical examples), we shall assume that:

$$\mathcal{F} = \mathcal{F}_\infty = \bigvee_{t \geq 0} \mathcal{F}_t.$$

Again, we will have to distinguish two cases: the case of arbitrary random times and honest times. Let $\rho$ be random time. We enlarge the initial filtration $(\mathcal{F}_t)$ with the process $(\rho \wedge t)_{t \geq 0}$, so that the new enlarged filtration $(\mathcal{F}^\rho_t)_{t \geq 0}$ is the smallest filtration (satisfying the usual assumptions) containing $(\mathcal{F}_t)$ and making $\rho$ a stopping time (i.e. $\mathcal{F}^\rho_t = \mathcal{K}^o_{t+}$, where $\mathcal{K}^o_t = \mathcal{F}_t \bigvee \sigma(\rho \wedge t)$). Sometimes it is more convenient to introduce the larger filtration

$$\mathcal{G}^\rho_t = \{A \in \mathcal{F}_\infty : \ \exists A_t \in \mathcal{F}_t,\ A \cap \{L > t\} = A_t \cap \{L > t\}\},$$

which coincides with $\mathcal{F}^\rho_t$ before $\rho$ and which is constant after $\rho$ and equal to $\mathcal{F}_\infty$ ([28], p. 186). In the case of an honest time $L$, one can show that in fact (see [41]):

$$\mathcal{F}^L_t = \{A \in \mathcal{F}_\infty : \ \exists A_t, B_t \in \mathcal{F}_t,\ A = (A_t \cap \{L > t\}) \cup (B_t \cap \{L \leq t\})\}.$$

**In the sequel, we shall only consider the filtrations $(\mathcal{G}^\rho_t)$ and $(\mathcal{F}^L_t)$: the first one when we study arbitrary random times and the second one when we consider the special case of honest times.**

*9.2.1. A description of predictable and optional processes in $\left(\mathcal{G}_t^\rho\right)$ and $\left(\mathcal{F}_t^L\right)$*

All the results we shall mention in what follows can be found in [43] (or in [42, 28]) and are particulary useful in mathematical finance ([30], [40]).

**Proposition 9.12.** *Let $\rho$ be an arbitrary random time. The following hold:*

1. *If $H$ is a $\left(\mathcal{G}_t^\rho\right)$ predictable process, then there exists a $\left(\mathcal{F}_t\right)$ predictable process $J$ such that*

$$H_t \mathbf{1}_{t \leq \rho} = J_t \mathbf{1}_{t \leq \rho}.$$

2. *If $T$ is a $\left(\mathcal{G}_t^\rho\right)$ stopping time, then there exists a $\left(\mathcal{F}_t\right)$ stopping time $S$ such that:*

$$T \wedge \rho = S \wedge \rho.$$

3. *Let $\xi \in L^1$. Then a càdlàg version of the martingale $\xi_t = \mathbb{E}\left[\xi | \mathcal{G}_t^\rho\right]$ is given by:*

$$\xi_t = \frac{1}{Z_t^\rho} \mathbf{1}_{t < \rho} \mathbb{E}\left[\xi \mathbf{1}_{t < \rho} | \mathcal{F}_t\right] + \xi \mathbf{1}_{t \geq \rho}.$$

**Proposition 9.13.** *Let $L$ be an honest time. The following hold:*

1. *$H$ is a $\left(\mathcal{F}_t^L\right)$ predictable process if and only if there exist two $\left(\mathcal{F}_t\right)$ predictable processes $J$ and $K$ such that*

$$H_t = J_t \mathbf{1}_{t \leq L} + K_t \mathbf{1}_{t > L}.$$

2. *Let $\xi \in L^1$. Then a càdlàg version of the martingale $\xi_t = \mathbb{E}\left[\xi | \mathcal{F}_t^L\right]$ is given by:*

$$\xi_t = \frac{1}{Z_t^L} \mathbb{E}\left[\xi \mathbf{1}_{t < L} | \mathcal{F}_t\right] \mathbf{1}_{t < L} + \frac{1}{1 - Z_t^L} \mathbb{E}\left[\xi \mathbf{1}_{t \geq L} | \mathcal{F}_t\right] \mathbf{1}_{t \geq L}.$$

3. *Every $\left(\mathcal{F}_t^L\right)$ optional process decomposes as*

$$H \mathbf{1}_{[0,L[} + J \mathbf{1}_{[L]} + K \mathbf{1}_{]L,\infty[},$$

*where $H$ and $K$ are $\left(\mathcal{F}_t\right)$ optional processes and where $J$ is a $\left(\mathcal{F}_t\right)$ progressively measurable process.*

**Proposition 9.14** ([43]). *Let $L$ be an honest time:*

1. *Let $X$ be a $\left(\mathcal{F}_t^L\right)$ local martingale stopped at $L$; then $X$ is also a $\left(\mathcal{G}_t^L\right)$ local martingale.*
2. *Let $Y$ be a $\left(\mathcal{G}_t^L\right)$ local martingale, stopped at $L$. If $Y_L$ is $\left(\mathcal{F}_L^L\right)$ measurable, then $Y$ is a $\left(\mathcal{F}_t^L\right)$ local martingale.*

Now we give a theorem which is very often used in applications:

**Theorem 9.15** ([43]). *Let $H$ be a bounded $\left(\mathcal{G}_t^\rho\right)$ predictable process. Then*

$$H_\rho \mathbf{1}_{\rho \leq t} - \int_0^{t \wedge \rho} \frac{H_s}{Z_{s-}^\rho} dA_s^\rho$$

*is a $\left(\mathcal{G}_t^\rho\right)$ martingale.*

**Remark 9.16.** If $L$ is an honest time, then

$$H_L \mathbf{1}_{L \le t} - \int_0^{t \wedge L} \frac{H_s}{Z_{s-}^L} dA_s^L$$

is a $\left( \mathcal{F}_t^L \right)$ martingale.

**Remark 9.17.** Let $\rho$ be a random time; taking $H \equiv 1$, we find that $\int_0^{t \wedge \rho} \frac{1}{Z_{s-}^{\rho}} dA_s^{\rho}$ is the $(\mathcal{G}_t^{\rho})$ dual predictable projection of $\mathbf{1}_{\rho \le t}$. When $\rho$ is a pseudo-stopping time that avoids $(\mathcal{F}_t)$ stopping times, we have from Theorem 8.32 that the $(\mathcal{G}_t^{\rho})$ dual predictable projection of $\mathbf{1}_{\rho \le t}$ is $\log \left( \frac{1}{Z_{t \wedge \rho}^{\rho}} \right)$.

Now, we shall study the properties $\rho$ as a stopping time in $(\mathcal{G}_t^{\rho})$.

**Proposition 9.18** ([43]).    *1. $\rho$ is a $(\mathcal{G}_t^{\rho})$ predictable stopping time if and only if it is a $(\mathcal{F}_t)$ predictable stopping time.*

   *2. Define*

$$\rho_1 = \rho \quad \text{on} \quad \Delta A_\rho^\rho = 0, \quad \rho_1 = \infty \quad \text{on} \quad \Delta A_\rho^\rho > 0,$$

     *and*

$$\rho_2 = \rho \quad \text{on} \quad \Delta A_\rho^\rho > 0, \quad \rho_2 = \infty \quad \text{on} \quad \Delta A_\rho^\rho = 0.$$

     *Then $\rho_1$ (resp. $\rho_2$) is the totally inaccessible part (resp. accessible part) of the $(\mathcal{G}_t^{\rho})$ stopping time $\rho$.*

   *3. If $L$ is an honest time, we can replace in the above $(\mathcal{G}_t^{\rho})$ with $\left( \mathcal{F}_t^L \right)$.*

*9.2.2. The decomposition formula before $\rho$*

In general, for an arbitrary random time, a local martingale is not a semimartingale in $(\mathcal{G}_t^{\rho})$. However, we have the following result:

**Theorem 9.19** (Jeulin-Yor [43]). *Every $(\mathcal{F}_t)$ local martingale $(M_t)$, stopped at $\rho$, is a $(\mathcal{G}_t^{\rho})$ semimartingale, with canonical decomposition:*

$$M_{t \wedge \rho} = \widetilde{M}_t + \int_0^{t \wedge \rho} \frac{d < M, \mu^\rho >_s}{Z_{s-}^{\rho}} \tag{9.8}$$

*where $\left( \widetilde{M}_t \right)$ is a $(\mathcal{G}_t^{\rho})$ local martingale.*

We shall now give two applications of this decomposition. The first one is a refinement of Theorem 8.32, which brings a new insight to pseudo-stopping times:

**Theorem 9.20.** *The following four properties are equivalent:*

   *1. $\rho$ is a $(\mathcal{F}_t)$ pseudo-stopping time, i.e (8.1) is satisfied;*
   *2. $\mu_t^\rho \equiv 1$, a.s*
   *3. $A_\infty^\rho \equiv 1$, a.s*

*4. every $(\mathcal{F}_t)$ local martingale $(M_t)$ satisfies*

$$(M_{t \wedge \rho})_{t \geq 0} \text{ is a local } (\mathcal{G}_t^\rho) \text{ martingale.}$$

*If, furthermore, all $(\mathcal{F}_t)$ martingales are continuous, then each of the preceding properties is equivalent to*

*5.*

$$(Z_t^\rho)_{t \geq 0} \text{ is a decreasing } (\mathcal{F}_t) \text{ predictable process}$$

*Proof.* $(1) \implies (2)$ For every square integrable $(\mathcal{F}_t)$ martingale $(M_t)$, we have

$$\mathbb{E}\left[M_\rho\right] = \mathbb{E}\left[\int_0^\infty M_s dA_s^\rho\right] = \mathbb{E}\left[M_\infty A_\infty^\rho\right] = \mathbb{E}\left[M_\infty \mu_\infty^\rho\right].$$

Since $\mathbb{E}M_\rho = \mathbb{E}M_0 = \mathbb{E}M_\infty$, we have

$$\mathbb{E}\left[M_\infty\right] = \mathbb{E}\left[M_\infty A_\infty^\rho\right] = \mathbb{E}\left[M_\infty \mu_\infty^\rho\right].$$

Consequently, $\mu_\infty^\rho \equiv 1$, *a.s*, hence $\mu_t^\rho \equiv 1$, *a.s* which is equivalent to: $A_\infty^\rho \equiv 1$, *a.s.* Hence, 2. and 3. are equivalent.

$(2) \implies (4)$. This is a consequence of the decomposition formula (9.8).

$(4) \implies (1)$ It suffices to consider any $\mathcal{H}^1$-martingale $(M_t)$, which, assuming (4), satisfies: $(M_{t \wedge \rho})_{t \geq 0}$ is a martingale in the enlarged filtration $(\mathcal{G}_t^\rho)$. Then as a consequence of the optional stopping theorem applied in $(\mathcal{G}_t^\rho)$ at time $\rho$, we get

$$\mathbb{E}\left[M_\rho\right] = \mathbb{E}\left[M_0\right],$$

hence $\rho$ is a pseudo-stopping time.

Finally, in the case where all $(\mathcal{F}_t)$ martingales are continuous, we show:

*a)* $(2) \Rightarrow (5)$ If $\rho$ is a pseudo-stopping time, then $Z_t^\rho$ decomposes as

$$Z_t^\rho = 1 - A_t^\rho.$$

As all $(\mathcal{F}_t)$ martingales are continuous, optional processes are in fact predictable, and so $(Z_t^\rho)$ is a predictable decreasing process.

*b)* $(5) \Rightarrow (2)$ Conversely, if $(Z_t^\rho)$ is a predictable decreasing process, then from the uniqueness in the Doob-Meyer decomposition, the martingale part $\mu_t^\rho$ is constant, i.e. $\mu_t^\rho \equiv 1$, *a.s.* Thus, 2 is satisfied. $\qquad \square$

Now, we apply the progressive enlargements techniques to the study of the Burkholder-Davis-Gundy inequalities. More precisely, what remains of the Burkholder-Davis-Gundy inequalities when stopping times $T$ are replaced by arbitrary random times $\rho$? The question of probabilistic inequalities at an arbitrary random time has been studied in depth by M. Yor (see [79], [81, 50] for details and references). For example, taking the special case of Brownian motion, it can easily be shown that there cannot exist a constant $C$ such that:

$$\mathbb{E}\left[|B_\rho|\right] \leq C \mathbb{E}\left[\sqrt{\rho}\right]$$

for any random time $\rho$. For if it were the case, we could take $\rho = \mathbf{1}_A$, for $A \in \mathcal{F}_\infty$, and we would obtain:

$$\mathbb{E}\left[|B_1| \, \mathbf{1}_A\right] \leq C\mathbb{E}\left[\mathbf{1}_A\right]$$

which is equivalent to: $|B_1| \leq C$, a.s., which is absurd. Hence it is not obvious that the "strict" BDG inequalities might hold for stopped local martingales at other random times than stopping times. However, we have the following positive result:

**Theorem 9.21** ([66]). *Let $p > 0$. There exist two universal constants $c_p$ and $C_p$ depending only on $p$, such that for any $(\mathcal{F}_t)$ local martingale $(M_t)$, with $M_0 = 0$, and any $(\mathcal{F}_t)$ pseudo-stopping time $\rho$ we have*

$$c_p\mathbb{E}\left[\left(< M >_\rho\right)^{\frac{p}{2}}\right] \leq \mathbb{E}\left[\left(M_\rho^*\right)^p\right] \leq C_p\mathbb{E}\left[\left(< M >_\rho\right)^{\frac{p}{2}}\right].$$

*Proof.* It suffices, with the previous Theorem, to notice that in the enlarged filtration $(\mathcal{G}_t^\rho)$, $(M_{t\wedge\rho})$ is a martingale and $\rho$ is a stopping time in this filtration; then, we apply the classical BDG inequalities. $\square$

**Remark 9.22.** The constants $c_p$ and $C_p$ are the same as those obtained for martingales in the classical framework; in particular the asymptotics are the same (see [17]).

**Remark 9.23.** It would be possible to show the above Theorem, just using the definition of pseudo-stopping times (as random times for which the optional stopping theorem holds); but the proof is much longer.

### 9.2.3. The decomposition formula for honest times

One of the remarkable features of honest time (discovered by Barlow [13]) is the fact that the pair of filtrations $\left(\mathcal{F}_t, \mathcal{F}_t^L\right)$ satisfies the $(H')$ hypothesis and every $(\mathcal{F}_t)$ local martingale $X$ is an $\left(\mathcal{F}_t^L\right)$ semimartingale. More precisely:

**Theorem 9.24.** *An $(\mathcal{F}_t)$ local martingale $(M_t)$, is a semimartingale in the larger filtration $\left(\mathcal{F}_t^L\right)$ and decomposes as:*

$$M_t = \widetilde{M}_t + \int_0^{t\wedge L} \frac{d\langle M, Z^L\rangle_s}{Z_{s_-}^L} - \int_L^t \frac{d\langle M, Z^L\rangle_s}{1 - Z_{s_-}^L}, \tag{9.9}$$

*where $\left(\widetilde{M}_t\right)_{t\geq 0}$ denotes a $\left(\left(\mathcal{F}_t^L\right), \mathbb{P}\right)$ local martingale.*

**Remark 9.25.** There are non honest times $\rho$ such that the pair $(\mathcal{F}_t, \mathcal{G}_t^\rho)$ satisfies the $(H')$ hypothesis: for example, the pseudo-stopping times of Proposition 8.38 enjoy this remarkable property (see [60]) (for other examples see [42]).

*Proof.* We shall give a proof under the conditions (**CA**), which are general enough for most of the applications. In this special case, it is an consequence of Theorem 9.9. Indeed, we saw in the course of the proof of Theorem 9.9 that (for ease of notations we drop the upper index $L$):

$$\mathbf{1}_{\{A_\infty > A_t\}} = \mathbf{1}_{\{L > t\}}, \text{ and } \mathbf{1}_{\{A_\infty = A_t\}} = \mathbf{1}_{\{L \leq t\}}.$$

Thus, by definition of $\mathcal{F}_t^L$, we have:

$$\mathcal{F}_t^L \subset \mathcal{F}_t^{\sigma(A_\infty)}.$$

Now, let $M$ be an $L^2$ bounded $(\mathcal{F}_t)$ martingale; the general case follows by localization. From Theorem 9.9

$$M_t = \widetilde{M}_t + \int_0^t \mathbf{1}_{\{L > s\}} \frac{d\langle M, \mu \rangle_s}{Z_s} - \int_0^t \mathbf{1}_{\{L \leq s\}} \frac{d\langle M, \mu \rangle_s}{1 - Z_s},$$

where $\left(\widetilde{M}_t\right)_{t \geq 0}$ denotes an $\left(\mathcal{F}_t^L\right)$ $L^2$ martingale. Thus, $\left(\widetilde{M}_t\right)$, which is equal to:

$$M_t - \left( \int_0^t \mathbf{1}_{\{L > s\}} \frac{d\langle M, \mu \rangle_s}{Z_s} - \int_0^t \mathbf{1}_{\{L \leq s\}} \frac{d\langle M, \mu \rangle_s}{1 - Z_s} \right),$$

is $\left(\mathcal{F}_t^L\right)$ adapted, and hence it is an $L^2$ bounded $\left(\mathcal{F}_t^L\right)$ martingale. $\qquad \square$

There are many applications of progressive enlargements of filtrations with honest times, but we do not have the place here to give them. At the end of this section, we shall give a list of applications and references. Nevertheless, we mention an extension of the BDG inequalities obtained by Yor:

**Proposition 9.26** (Yor [81], p. 57). *Assume that $(\mathcal{F}_t)$ is the filtration of a standard Brownian Motion and let $L$ be an honest time. Then we have:*

$$\mathbb{E}\left[|B_L|\right] \leq C \mathbb{E}\left[\Phi_L \sqrt{L}\right],$$

*with*

$$\Phi_L = \left(1 + \log \frac{1}{I_L}\right)^{1/2}, \quad \text{where } I_L = \inf_{u < L} Z_u^L,$$

*and $C$ a universal constant.*

**Remark 9.27.** If $L$ is a stopping time, then $\Phi_L = 1$. Furthermore, for any continuous increasing function $f : \mathbb{R}_+ \to \mathbb{R}_+$, we have:

$$\mathbb{E}\left[f\left(\Phi_L\right)\right] \leq \mathbb{E}\left[f\left(V\right)\right],$$

with $V = (1 + \mathbf{e})^{1/2}$, where $\mathbf{e}$ is an exponential random variable with parameter 1.

**Remark 9.28.** It is also possible to prove that the BDG inequalities never hold for honest times under (**CA**) (see [66]).

### *9.3. The (H) hypothesis*

In this paragraph, we shall briefly mention the $(H)$ hypothesis, which is very widely used in the models of default times in mathematical finance. Let $(\Omega, \mathcal{F}, \mathbb{P})$ be a probability space satisfying the usual assumptions. Let $(\mathcal{F}_t)$ and $(\mathcal{G}_t)$ be two sub-filtrations of $\mathcal{F}$, with

$$\mathcal{F}_t \subset \mathcal{G}_t.$$

**Theorem 9.29** (Brémaud and Yor [22])**.** *The following are equivalent:*

1. *$(H)$: every $(\mathcal{F}_t)$ martingale is a $(\mathcal{G}_t)$ martingale;*
2. *for all $t \geq 0$, the sigma fields $\mathcal{G}_t$ and $\mathcal{F}_\infty$ are independent conditionally on $\mathcal{F}_t$.*

**Remark 9.30.** We shall also say that $(\mathcal{F}_t)$ is immersed in $(\mathcal{G}_t)$.

Now let us consider the $(\mathbf{H})$ hypothesis in the framework of the progressive enlargement of some filtration $(\mathcal{F}_t)$ with a random time $\rho$. This problem was studied by Dellacherie and Meyer [25]. It is equivalent to one of the following hypothesis (see [30] for more references):

1. $\forall t$, the $\sigma$-algebras $\mathcal{F}_\infty$ and $\mathcal{F}_t^\rho$ are conditionally independent given $\mathcal{F}_t$.
2. For all bounded $\mathcal{F}_\infty$ measurable random variables $\mathbf{F}$ and all bounded $\mathcal{F}_t^\rho$ measurable random variables $\mathbf{G}_t$, we have

$$\mathbb{E}\left[\mathbf{F}\mathbf{G}_t \mid \mathcal{F}_t\right] = \mathbb{E}\left[\mathbf{F} \mid \mathcal{F}_t\right] \mathbb{E}\left[\mathbf{G}_t \mid \mathcal{F}_t\right].$$

3. For all bounded $\mathcal{F}_t^\rho$ measurable random variables $\mathbf{G}_t$:

$$\mathbb{E}\left[\mathbf{G}_t \mid \mathcal{F}_\infty\right] = \mathbb{E}\left[\mathbf{G}_t \mid \mathcal{F}_t\right].$$

4. For all bounded $\mathcal{F}_\infty$ measurable random variables $\mathbf{F}$,

$$\mathbb{E}\left[\mathbf{F} \mid \mathcal{F}_t^\rho\right] = \mathbb{E}\left[\mathbf{F} \mid \mathcal{F}_t\right].$$

5. For all $s \leq t$,

$$\mathbb{P}\left[\rho \leq s \mid \mathcal{F}_t\right] = \mathbb{P}\left[\rho \leq s \mid \mathcal{F}_\infty\right].$$

Now we come back to the general situation described in Theorem 9.29. We assume that the hypothesis $(H)$ holds. What happens when we make an equivalent change of probability measure?

**Proposition 9.31** ([44])**.** *Let $\mathbb{Q}$ be a probability measure which is equivalent to $\mathbb{P}$ (on $\mathcal{F}$). Then every $(\mathcal{F}_\bullet, \mathbb{Q})$ semimartingale is a $(\mathcal{G}_\bullet, \mathbb{Q})$ semimartingale.*

Now, define:

$$\frac{d\mathbb{Q}}{d\mathbb{P}}\Big|_{\mathcal{F}_t} = R_t; \quad \frac{d\mathbb{Q}}{d\mathbb{P}}\Big|_{\mathcal{G}_t} = R_t'.$$

Jeulin and Yor [44] prove the following facts: if $Y = \dfrac{d\mathbb{Q}}{d\mathbb{P}}$, then the hypothesis $(H)$ holds under $\mathbb{Q}$ if and only if:

$$\forall X \geq 0,\ X \in \mathcal{F}_\infty, \quad \frac{\mathbb{E}_{\mathbf{P}}\left[XY|\mathcal{G}_t\right]}{R_t'} = \frac{\mathbb{E}_{\mathbf{P}}\left[XY|\mathcal{F}_t\right]}{R_t}.$$

In particular, when $\dfrac{d\mathbb{Q}}{d\mathbb{P}}$ is $\mathcal{F}_\infty$ measurable, $R_t = R'_t$ and the hypothesis $(H)$ holds under $\mathbb{Q}$.

Now let us give a decomposition formula:

**Theorem 9.32** (Jeulin-Yor [44])**.** *If $(X_t)$ is a $(\mathcal{F}_\bullet, \mathbb{Q})$ local martingale, then the stochastic process:*

$$I_X(t) = X_t + \int_0^t \frac{R'_{s-}}{R'_s}\left(\frac{1}{R_{s-}}d[X,R]_s - \frac{1}{R'_{s-}}d[X,R']_s\right)$$

*is a $(\mathcal{G}_\bullet, \mathbb{Q})$ local martingale.*

### 9.4. Concluding remarks on enlargements of filtrations

The theory of enlargements of filtrations has many applications and it is of course impossible to expose them in such an essay. I shall here simply mention some of its applications with references. Of course, the following list is far from being complete.

1. The theory of enlargements of filtrations is very efficient for proving path decompositions results (which are very often difficult to establish); the reader can refer to [42],[81],[60] or [63].
2. It can also be useful to obtain canonical decompositions of some bridges (just as we did for the Brownian Motion); see for example [81].
3. The theory of enlargements of filtrations is also very important in the study of the stopping theorems with stopping times replaced with arbitrary random times and in the study of zeros of continuous martingales: [12], [6],[8], [60],[81].
4. It is an important tool in mathematical finance, in the modeling of default times and in insider trading models: [30], [40], [33].
5. It is also useful to obtain general probabilistic inequalities: see [79], [81, 50].
6. Enlargements of filtrations sometimes bring new insight on some already known results such as Pitman's theorem, or Perkins's theorem on Brownian local times: [45].